\documentclass[a4paper,10pt]{article}
\usepackage[utf8x]{inputenc}
\usepackage{graphicx}
\usepackage{epsfig}
\usepackage{amssymb}
\usepackage{amsmath}
\usepackage{amsfonts}
\usepackage{moreverb}
\usepackage{rotating}
\usepackage{multirow}

\usepackage[colorlinks,bookmarksopen,bookmarksnumbered,citecolor=red,urlcolor=red]{hyperref}
\hypersetup{pdftitle={SWASHES: a compilation of Shallow Water Analytic Solutions for Hydraulic and Environmental Studies},
bookmarks=true,
pdftoolbar=true,
pdfmenubar=true,
pdfauthor={O. Delestre, C. Lucas, P.-A. Ksinant, F. Darboux, C. Laguerre, T.N.T. Vo, F. James, S. Cordier},
pdfsubject={applied mathematics, hydraulics},
pdfcreator={Delestre Olivier},
pdfproducer={Delestre Olivier},
pdfkeywords={Shallow-Water equation} {Saint-Venant} {Analytic solutions} {benchmarking} {SWASHES}
 {source terms} {friction} {rain} {Manning} {Chezy} {Strickler} {Darcy-Weisbach} {Steady-state flow}
{transitory flow} {Validation of numerical methods} {wet/dry transition}}

\usepackage{multirow}

\newtheorem{remark}{Remark}
\def\Fr{{\rm Fr}}

\title{SWASHES: a compilation of Shallow Water Analytic Solutions for Hydraulic and Environmental Studies}
\author{{{O. Delestre}\footnote{Corresponding author: Delestre@unice.fr, presently at: Laboratoire de Math\'ematiques
 J.A. Dieudonn\'e -- Polytech Nice-Sophia , Universit\'e de Nice -- Sophia Antipolis, Parc Valrose, 06108 Nice cedex 02,
 France}$^{\; ,\ddag,}$}\footnote{Institut Jean le Rond d'Alembert, CNRS \& UPMC Universit\'e Paris 06, UMR 7190,
 4 place Jussieu, Bo\^ite 162, F-75005 Paris, France}, C. Lucas\footnote{MAPMO UMR CNRS 7349, Universit\'e d'Orl\'eans,
 UFR Sciences, B\^atiment de math\'ematiques,
 B.P. 6759 -- F-45067 Orl\'eans cedex 2, France}\;,
 {P.-A. Ksinant$^{\ddag,}$}\footnote{INRA, UR 0272 Science du sol,
 Centre de recherche d'Orl\'eans, CS 40001 Ardon, F-45075 Orl\'eans cedex 2, France}\,, F. Darboux$^\S$,\\
 C. Laguerre$^{\ddag}$,  {T.N.T. Vo$^{\ddag,}$}\footnote{Presently at: Department of Applied Mathematics,
 National University of Ireland, Galway, Republic of Ireland}\,, F. James$^{\ddag}$ and S. Cordier$^{\ddag}$}

\begin{document}

\maketitle

\begin{abstract}
Numerous codes are being developed to solve Shallow Water equations.
Because there are used in hydraulic and environmental studies, their capability to simulate properly flow dynamics is critical to guarantee infrastructure and human safety.
While validating these codes is an important issue, code validations are currently restricted because analytic solutions to the Shallow Water equations are rare and have been published on an individual basis over a period of more than five decades.
This article aims at making analytic solutions to the Shallow Water equations easily available to code developers and users.
It compiles a significant number of analytic solutions to the Shallow Water equations
that are currently scattered through the literature of various 
scientific disciplines.
The analytic solutions are described in a unified formalism to make a consistent set of test cases.
These analytic solutions encompass a wide variety of
flow conditions (supercritical, subcritical, shock, etc.),
in 1 or 2 space dimensions, with or without rain and
soil friction, for transitory flow or steady state.
The corresponding source codes
are made available to the community (\url{http://www.univ-orleans.fr/mapmo/soft/SWASHES}), so that 
users of Shallow Water-based models can easily find an adaptable benchmark library to
validate their numerical methods.

\end{abstract}



\section{Introduction}

Shallow-Water equations have been proposed by Adh\'emar Barr\'e de Saint-Venant in 1871 
to model flows in a channel \cite{saintvenant71}. 
Nowadays,  they are widely used to model flows in various contexts, such as: overland flow \cite{Esteves00, Tatard08},
rivers \cite{Goutal02, Burguete04}, flooding \cite{Caleffi03, Delestre09},
dam breaks \cite{Alcrudo99, Valiani02}, 
nearshore \cite{Borthwick01, Marche05}, 
tsunami \cite{George06, Kim07, Popinet11}.
These equations consist in a nonlinear system of partial differential equations (PDE-s), more precisely
conservation laws describing the evolution of the height and mean velocity of the fluid.

In real situations (realistic geometry, sharp spatial or temporal variations of the parameters in the model, etc.),
there is no hope to solve explicitly this system of PDE-s, \emph{i.e.}
to produce analytic formul\ae \ for the solutions.
It is therefore necessary to develop specific numerical methods to compute approximate solutions of such PDE-s, 
see for instance \cite{Toro97, LeVeque02, Bouchut04}.
Implementation of any of such methods raises the question of the validation of the code.  

Validation is an essential step to check if a model 
(that is the equations, the numerical methods and their implementation)
suitably describes the considered phenomena.
There exists at least three complementary types of numerical tests to
ensure a numerical code is relevant for the considered systems of equations.
First, one can produce convergence or stability results (\emph{e.g.} by refining the mesh). 
This validates only the numerical method and its implementation.
Second, approximate solutions can be matched to analytic solutions available for some
simplified or specific cases. 
Finally, numerical results can be compared with experimental data, produced indoor
or outdoor. This step should be done after the previous two; it is the most difficult one
and must be validated by a specialist of the domain. 
This paper focuses on the second approach.

Analytic solutions seem underused in the validation of numerical codes, possibly for the following reasons.
First, each analytic solution has a limited scope in terms of flow conditions.
Second, they are currently scattered through the literature and, thus, are difficult to find.
However, there exists a significant number of published analytic solutions that encompasses a wide range of flow conditions.
Hence, this gives a large potential to analytic solutions for validatation of numerical codes.

This work aims at overcoming these issues, on the one hand by gathering a significant set of analytic solutions,
on the other hand by providing the corresponding source codes. 
The present paper describes the analytic solutions together with some comments about their interest and use.
The source codes are made freely available to the community through the SWASHES 
(Shallow-Water Analytic Solutions for Hydraulic and Environmental Studies) software. 
The SWASHES library does not pretend to list all available analytic solutions.
On the contrary, it is open for extension and we take here the opportunity to ask users to contribute to the project by 
sending other analytic solutions together with the dedicated code.

The paper is organized as follows: in Section \ref{sec:eq}, we briefly present the notations we use and the
main properties of Shallow-Water equations. 
In the next two sections, we briefly outline each analytic solution.
A short description of the SWASHES software can be find in Section \ref{sec:SWASHES}.
The final section is an illustration using the results of the Shallow-Water code 
developed by our team for a subset of analytic solutions.

\section{Equations, notations and properties }\label{sec:eq}

First we describe the rather general settings of viscous Shallow-Water equations
in two space dimensions, with topography, rain, infiltration and soil friction. 
In the second paragraph, we give the simplified system arising in one space dimension and recall several
classical properties of the equations.

\subsection{General settings}

\begin{figure}
\begin{center}
\includegraphics[width = 0.5 \textwidth]{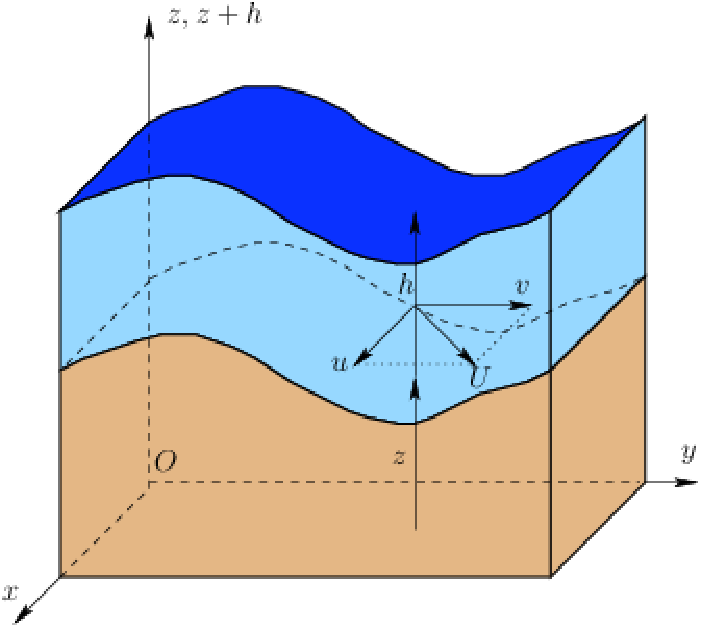}
\caption{Notations for 2D Shallow-Water equations}
\label{fignot2d}
\end{center}
\end{figure}

The unknowns of the equations are the water height ($h(t,x,y)$ [L]) and $u(t,x,y)$, $v(t,x,y)$ 
the horizontal components of the vertically averaged velocity [L/T] (Figure~\ref{fignot2d}). 

The equations take the following form of balance laws, 
where $g=9.81$ m/s$^2$ is the gravity constant:
\begin{equation}\label{SaintVenant2D}
\left\{
\begin{array}{l}
\partial_t h +\partial_x\left(hu\right)+\partial_y\left(hv\right)=R-I\\
\\
\partial_t\left(hu\right)+\partial_x\left(hu^2+\dfrac{gh^2}{2}\right)+\partial_y\left(huv\right)
  = gh({S_0}_x-{S_f}_x) + \mu {S_d}_x\\
\\
\partial_t\left(hv\right)+\partial_x\left(huv\right)+\partial_y\left(hv^2+\dfrac{gh^2}{2}\right)
  = gh({S_0}_y-{S_f}_y) + \mu {S_d}_y
\end{array}
\right.
\end{equation}
The first equation is actually a mass balance. The fluid density can be replaced by the height because of 
incompressibility. The other two equations are momentum balances, and involve forces such as gravity and friction.
We give now a short description of all the terms involved, recalling the physical dimensions.

\begin{itemize}
\item $z$ is the topography [L], since we consider no erosion 
here, it is a fixed function of space,
$z(x,y)$, and we classically denote by $S_{0x}$ (resp. $S_{0y}$) the opposite of the slope in the $x$ (resp. $y$) direction,
$S_{0x}=-\partial_x z(x,y)$ (resp. $S_{0y}=-\partial_y z(x,y)$);
\item $R\geq 0$ is the rain intensity [L/T], it is a given function $R(t,x,y)\ge 0$. In this paper, it is considered uniform in  space;
\item $I$ is the infiltration rate [L/T], mentioned for the sake of completeness. It is given by another
 model (such as Green-Ampt, Richards, etc.) and is not taken into account in the following;
\item $S_f=\left({S_f}_x,{S_f}_y\right)$ is the friction force. The friction law $S_f$ may take several forms, depending on both soil and flow properties.
In the formul\ae\ below, $U$ is the velocity vector $U=(u,v)$ with $|U| = \sqrt {u^2 + v^2}$ and $Q$ is the discharge $Q=(hu,hv)$.
In hydrological models, two families of friction laws are encountered, based on empirical considerations.
On the one hand, we have the family of Manning-Strickler's friction laws
\begin{equation*}
 S_f=C_f\dfrac{U|U|}{h^{4/3}}=C_f\dfrac{Q|Q|}{h^{10/3}}
\end{equation*}
 $C_f=n^2$, where $n$ is the Manning's coefficient [L\textsuperscript{-1/3}T].\\
On the other hand, the laws of Darcy-Weisbach's and Ch\'ezy's family writes
\begin{equation*}
 S_f=C_f\dfrac{U|U|}{h}=C_f\dfrac{Q|Q|}{h^{3}}.
\end{equation*}
With $C_f=f/(8g)$, $f$ dimensionless coefficient, (resp. $C_f=1/{C^2}$, $C$ [L\textsuperscript{1/2}/T]) 
we get the Darcy-Weisbach's (resp. Ch\'ezy's) friction law.
Notice that the friction may depend on the space variable, especially on large parcels. In the sequel this will
not be the case.
\item finally, $\mu S_d = \left(\mu {S_d}_x, \mu {S_d}_y \right)$ is the viscous term with $\mu\ge 0$ the viscosity of the fluid [L\textsuperscript{2}/T].
\end{itemize}

\subsection{Properties}

In this section, we recall several properties of the Shallow-Water model that are useful in the flow description.
For the sake of simplicity, we place ourselves in the one-dimensional case, extensions to the general setting being straightforward.
The two-dimensional equations \eqref{SaintVenant2D} rewrite
\begin{equation}\label{SaintVenant1D}
 \left\{\begin{array}{l}
        \partial_t h+\partial_x (hu) = R-I\\
  \partial_t (hu) +\partial_x \left(hu^2+\dfrac{gh^2}{2}\right)=gh(S_{0x}-S_f)+\mu \partial_x\left(h\partial_x u\right)
       \end{array}\right.
\end{equation}

The left-hand side of this system is the transport operator, corresponding to the flow of an ideal fluid in a flat channel, without friction, 
rain or infiltration. This is actually the model introduced by Saint-Venant in \cite{saintvenant71}, and it contains several important properties of the flow. 
In order to emphasize these properties, we first rewrite the one-dimensional equations using vectors form:
\begin{equation}\label{nonconservative}
\partial_t W+\partial_x F(W) = 0, \quad\mbox{where } W=\left(\begin{array}{c}h\\hu\end{array}\right),\quad 
  F(W)=\left(\begin{array}{c}hu\\hu^2+\dfrac{gh^2}{2}\end{array}\right),
\end{equation}
with $F(W)$ the flux of the equation. The transport is more clearly evidenced in the following nonconservative form,
where $A(W)=F'(W)$ is the matrix of transport coefficients:
  \begin{equation}\label{MatrixTransp}
\partial_t W+A(W)\partial_x W = 0, \quad A(W)=F'(W)=\left(\begin{array}{cc}0&1\\-u^2+gh&2u\end{array}\right).
  \end{equation}
More precisely, when $h>0$, the matrix $A(W)$ turns out to be diagonalizable, with eigenvalues
  $$
\lambda_1(W) = u-\sqrt{gh} < u+\sqrt{gh} = \lambda_2(W).
  $$
This important property is called strict hyperbolicity (see for instance \cite{Godlewski96} and references therein for more complete
information). The eigenvalues are indeed velocities, namely the ones of surface waves on the fluid, which are basic characteristics of
the flow. Notice here that the eigenvalues
coincide if $h=0$~m, that is for dry zones. In that case, the system is no longer hyperbolic, and this induces difficulties at both theoretical and
numerical levels. Designing numerical schemes that preserve positivity for $h$ is very important in this context.

From these formul\ae\ we recover a useful classification of flows, based on the relative values of the velocities of the fluid, $u$, and of the 
waves, $\sqrt{gh}$. Indeed if  $|u| < \sqrt{gh}$ the characteristic velocities have opposite signs, and information propagate upward as
well as downward the flow, which is then said subcritical or fluvial. On the other hand, when $|u| > \sqrt{gh}$,
the flow is supercritical, or torrential, all the information go downwards. 
A transcritical regime exists when some parts of a flow are subcritical, other supercritical.

Since we have two unknowns $h$ and $u$ (or equivalently $h$ and $q=hu$), a subcritical flow is therefore determined by one upstream and one downstream value, whereas a supercritical flow is completely determined by the two upstream values. 
Thus for numerical simulations, we have to impose one variable for subcritical inflow/outflow. We impose both variables for supercritical
inflow and for supercritical outflow, free boundary conditions are considered (see for example \cite{Bristeau01}).

In this context, two quantities are useful. The first one is a dimensionless parameter called the Froude number
  \begin{equation}\label{DefFroude}
\Fr = \dfrac{|u|}{\sqrt{gh}}.
  \end{equation}
It is the analogue of the Mach number in gas dynamics, and the flow is subcritical (resp. supercritical) if $\Fr<1$ (resp. $\Fr>1$).
The other important quantity is the so-called critical height $h_c$ which writes
  \begin{equation}\label{DefCritHeight}
h_c = \left(\dfrac{q}{\sqrt{g}}\right)^{2/3},
  \end{equation}
for a given discharge $q=hu$.
It is a very readable criterion for criticality: the flow is  subcritical (resp. supercritical) if $h>h_c$ (resp. $h<h_c$).

When additional terms are present, other properties have to be considered,
for instance the occurrence of steady states (or equilibrium) solutions. 
These specific flows are defined and discussed in Section~\ref{sec:steady-state}.

\label{MoreCells}

\section{Steady state solutions}\label{sec:steady-state}

In this section, we focus on a family of steady state solutions, that is solutions that satisfy:
\begin{equation*}
 \partial_t h=\partial_t u=0.
\end{equation*}
Replacing this relation in the one-dimensional Shallow-Water equations \eqref{SaintVenant1D},
the mass equation gives $\partial_x (hu)=R$ or $hu=q=Rx+q_0$, where $q_0=q(t,x=0)$.
Similarly the momentum equation writes
\begin{equation*}
 \partial_x\left(\dfrac{q^2}{h}+ \dfrac{gh^2}{2}\right)=-gh\partial_x z-ghS_f(h,q)+\mu \partial_x\left(h \partial_x \dfrac{q}{h}\right).
\end{equation*}
Thus for $h\neq 0$, we have the following system
\begin{equation}
 \left\{\begin{array}{l}
         q=Rx+q_0,\\
  \partial_x z=\dfrac{1}{gh}\left(\dfrac{q^2}{h^2}-gh\right)\partial_x h-S_f(h,q)
+\dfrac{\mu}{gh}\partial_x\left(h\partial_x \dfrac{q}{h}\right).
        \end{array}\right.\label{eq:general-steady-SV}
\end{equation}
System \eqref{eq:general-steady-SV} is the key point of the following series of analytic solutions. For these solutions,
the strategy consists in choosing either a topography and getting the associated water height or a water height and deducing
the associated topography. \\

Since \cite{Bermudez94}, it is well known that the source term treatment is a crucial point in preserving steady states.
With the following steady states solutions, one can check
if the steady state at rest and dynamic steady states are satisfied 
by the considered schemes using various flow conditions (fluvial, torrential,
transcritical, with shock, etc.). Moreover, 
the variety of inflow and outflow configurations (flat bottom/varying topography,
with/without friction, etc.)  gives a validation of boundary conditions treatment.
One must note that, as different source terms (topography, friction, rain and diffusion)
are taken into account, these solutions can also validate the source terms treatment.\\
The last remark deals with initial conditions: if initial conditions are taken equal to the solution at the steady state,
one can only conclude on the ability of the numerical scheme to \emph{preserve} steady states. In order to prove the
capacity to \emph{catch} these states, initial conditions should be different from the steady state. This is the reason
why the initial conditions, as well as the boundary conditions, are described in each case.  

Table~\ref{AnaSolTable-Steady} lists all steady-state solutions available in SWASHES and outlines their main features.

\begin{table}[htbp]
{\hspace*{-4cm}
\begin{tiny}
\setlength{\tabcolsep}{1mm}
\begin{tabular}{|c|c|c|c|c|c|c|c|c|c|c|c|c|}  \cline{1-12}

\multicolumn{4}{|l|}{\emph{Steady-state solutions}}&
\multicolumn{4}{c|}{Flow criticality}&
\multicolumn{4}{c|}{Friction} &
 \multicolumn{1}{c}{} \\[0.1cm] \hline
\multicolumn{1}{|c|}{Type} &\multicolumn{1}{|c|}{Description}&\multicolumn{1}{|c|}{\textsection} & \multicolumn{1}{c|}{Reference} &
                     Sub. & Sup. & Sub.$\rightarrow$Sup.& Jump & Man.& D.-W. & Other & Null & Comments \\\hline

\multirow{6}{*}{Bumps}
 
&
Lake at rest with immersed bump
& \ref{Bump-PuddleImmersed}
&  \cite{Delestre10b}
&     &     &                    &  & & & & X &
Hydrostatic equilibria\\

\cline{2-13}

&
Lake at rest with emerged bump
&\ref{Bump-PuddleEmerged}
&  \cite{Delestre10b}           
&     &     &                    &  & & & & X &
\begin{tabular}{c}Hydrostatic equilibria and\\wet-dry transition\end{tabular} \\

\cline{2-13}

& Subcritical flow
&\ref{Bump-Subcritical}
& \cite{Goutal97}
&  X  &     &                    &  & & & & X &
Initially steady state at rest\\

\cline{2-13}

&
Transcritical flow without shock
&\ref{Bump-Transcritical-NoShock}
& \cite{Goutal97}
&     &     &  X                 &  & & & & X &
Initially steady state at rest\\

\cline{2-13}

&
Transcritical flow with shock
&\ref{Bump-Transcritical-WithShock-NoFriction}
& \cite{Goutal97}
&     &     &                    & X  & & & & X &
Initially steady state at rest\\

\hline

\multirow{35}{*}{
\begin{tabular}{c}
Flumes\\
(MacDonald's\\
based)
\end{tabular}
}

&\begin{tabular}{c}Long channel\\with subcritical flow\end{tabular}
&\ref{McDo-1D-1000m-Subcritical}
& \cite{Vo08}
&  X  &     &    &    &   X  &  X  & & &  Initially dry channel. 1000~m long\\  

\cline{2-13}

&\begin{tabular}{c}Long channel\\with supercritical flow\end{tabular}
&\ref{McDo-1D-1000m-Supercritical}
& \cite{Delestre10b}
&    &  X   &    &    &   X  &  X  & & & Initially dry channel.  1000~m long\\  

\cline{2-13}

&\begin{tabular}{c}Long channel\\with sub- to super-critical flow\end{tabular}
&\ref{McDo-1D-1000m-SubSupercritical}
& \cite{Vo08}
&    &     &  X  &    &  X   &  X   & & & Initially dry channel.  1000~m long\\  

\cline{2-13}

&\begin{tabular}{c}Long channel\\with super- to sub-critical flow\end{tabular}
&\ref{McDo-1D-1000m-SuperSubcritical}
& \cite{Vo08}
&    &     &    &  X  &   X  &  X  & & &  Initially dry channel.  1000~m long\\  

\cline{2-13}

&\begin{tabular}{c}Short channel\\with smooth transition and shock\end{tabular}
&\ref{McDo-1D-100m-Subcritical}
& \cite{Vo08}
&   &     & X   &  X  &   X  &   &  & &  
\begin{tabular}{c}At t=0, lake downstream.\\100~m long\end{tabular}\\  

\cline{2-13}

&\begin{tabular}{c}Short channel\\with supercritical flow\end{tabular}
&\ref{McDo-1D-100m-Supercritical}
&  \cite{Delestre10b}
&    &  X   &    &    &  X   &   &  & & Initially dry. 100~m long\\  

\cline{2-13}

&\begin{tabular}{c}Short channel\\with sub- to super-critical flow\end{tabular}
&\ref{McDo-1D-100m-SubSupercritical}
& \cite{Vo08}
&    &     & X  &    &  X   &    & & &  
\begin{tabular}{c}At t=0, lake downstream.\\100~m long\end{tabular}\\  

\cline{2-13}

&\begin{tabular}{c}Very long, undulating and periodic channel\\with subcritical flow\end{tabular}
&\ref{McDo-1D-5000m-Periodic-Subcritical}
& \cite{Vo08}&  X  &     &    &    &  X   &   &  & &
\begin{tabular}{c}At t=0, lake downstream.\\5000~m long\end{tabular}\\  

\cline{2-13}

&\begin{tabular}{c}Rain on a long channel\\with subcritical flow\end{tabular}
&\ref{McDo-1Drain-Subcritical}
& \cite{Vo08}
&  X  &     &    &    &  X   &  X  & & &  Initially dry. 1000~m long\\  

\cline{2-13}

&\begin{tabular}{c}Rain on a long channel\\with supercritical flow\end{tabular}
&\ref{McDo-1Drain-Supercritical}
& \cite{Vo08}
&    &  X   &    &    &  X   &  X  & & &  Initially dry. 1000~m long \\  

\cline{2-13}

&\begin{tabular}{c}Long channel\\with subcritical flow and diffusion\end{tabular}
&\ref{McDo-1Ddiffusion-Subcritical}
& \cite{Delestre10}
&  X  &     &    &    &    &     & X & & Initially dry. 1000~m long  \\  

\cline{2-13}

&\begin{tabular}{c}Long channel\\with supercritical flow and diffusion\end{tabular}
&\ref{McDo-1Ddiffusion-Supercritical}
& \cite{Delestre10}
&    &  X   &    &    &     &    & X & &  Initially dry. 1000~m long \\

\cline{2-13}

&\begin{tabular}{c}Pseudo-2D short channel\\with subcritical flow\end{tabular}
&\ref{McDo-pseudo2D-Subcritical}
& \cite{MacDonald96}
&  X  &     &    &    &   X  &   &  & &
\begin{tabular}{c}Rectangular cross section.\\Initially dry. 200~m long\end{tabular}\\   

\cline{2-13}

&\begin{tabular}{c}Pseudo-2D short channel\\with supercritical flow\end{tabular}
&\ref{McDo-pseudo2D-Supercritical}
& \cite{MacDonald96}
 &    &  X   &    &    &    X &   &  & &
\begin{tabular}{c}Rectangular cross section.\\Initially wet. 200~m long\end{tabular}\\   

\cline{2-13}

&\begin{tabular}{c}Pseudo-2D short channel\\with smooth transition\end{tabular}
&\ref{McDo-pseudo2D-Smooth}
& \cite{MacDonald96}
&    &     &  X  &    &   X  &   &  & &
\begin{tabular}{c}Rectangular cross section.\\Initially partly-wet. 200~m long\end{tabular}\\   

\cline{2-13}

&\begin{tabular}{c}Pseudo-2D short channel\\with shock\end{tabular}
&\ref{McDo-pseudo2D-Jump}
& \cite{MacDonald96}
&    &     &    &  X  &   X  &   &  & &
\begin{tabular}{c}Rectangular cross section.\\Initially dry. 200~m long\end{tabular}\\   

\cline{2-13}

&\begin{tabular}{c}Pseudo-2D long channel\\with subcritical flow\end{tabular}
&\ref{McDo-pseudo2D-Subcritical2}
& \cite{MacDonald96}
&  X  &     &    &    &  X   &   &  & &
\begin{tabular}{c}Isoscele trapezoidal cross section.\\Initially dry. 400~m long\end{tabular}\\   

\cline{2-13}

&\begin{tabular}{c}Pseudo-2D long channel\\with smooth transition and shock\end{tabular}
&\ref{McDo-pseudo2D-SmoothJump}
& \cite{MacDonald96}
 &    &    &  X  & X   &  X   &   &  & &  
\begin{tabular}{c}Isoscele trapezoidal cross section.\\Initially dry. 400~m long\end{tabular}\\   

\hline
\multicolumn{13}{l}{\emph{Slope is always variable. Sub.: Subcritical;  Sup.: Supercritical; Man.: Manning; D.-W.: Darcy-Weisbach}}\\

\end{tabular}
\end{tiny}

{\caption{Analytic solutions for shallow flow equations and their main features --- Steady-state cases
\label{AnaSolTable-Steady} }}

}
\end{table}


\subsection{Bumps}

Here we present a series of steady state cases proposed in \cite[p.14-17]{Goutal97} based on an idea introduced in
 \cite{Houghton68}, with a flat topography at the boundaries, no rain, no friction and no diffusion ($R=0$~m/s, $S_f=0$
 and $\mu=0~\text{m}^2/\text{s}$). Thus system \eqref{eq:general-steady-SV} reduces to
\begin{equation*}
 \left\{\begin{array}{l}
         q=q_0,\\
    \partial_x z=\dfrac{1}{gh}\left(\dfrac{q^2}{h^2}-gh\right)\partial_x h.
\end{array}\right.\label{eq:steady-SV-bump}
\end{equation*}
In the case of a regular solution, we get the Bernoulli relation
\begin{equation}
 \dfrac{q_0^2}{2gh^2(x)}+h(x)+z(x)=Cst\label{eq:Bernoulli}
\end{equation}
which gives us the link between the topography and the water height.

\noindent
Initial conditions satisfy the hydrostatic equilibrium
\begin{equation}
h+z=Cst\;\text{and}\;q=0~\text{m}^2/\text{s}.\label{hydeq}
\end{equation}
These solutions test the preservation of steady states and the boundary conditions treatment.

\medskip

In the following cases, we choose a domain of length $L = $~25~m with a topography given by:
\begin{equation*}
z(x)=\left\{\begin{array}{ll}
0.2-0.05(x-10)^2&\text{if}\;8\;\text{m}<x<12\;\text{m},\\
0&\text{else}. \end{array}\right.
 \end{equation*}


\subsubsection{Lake at rest with an immersed bump}
\label{Bump-PuddleImmersed}
In the case of a lake at rest with an immersed bump, the water height is such that the topography is totally immersed \cite{Delestre10b}.
In such a configuration, starting from the steady state, the velocity must be null and the water surface should stay flat.\\

In SWASHES we have the following initial conditions: \\
\[h+z=0.5\;\text{m and}\;q=0\;\text{m}^2\text{/s}\]
and the boundary conditions
\begin{equation*}
\left\{\begin{array}{l}
h=0.5\;\text{m},\\
q=0\;\text{m}^2\text{/s}.\end{array}\right.
 \end{equation*}
 

\subsubsection{Lake at rest with an emerged bump}
\label{Bump-PuddleEmerged}
The case of a lake at rest with an emerged bump is the same as in the previous section except that the water height is smaller
in order to have emergence of some parts of the topography \cite{Delestre10b}.   
Here again, we initialize the solution at steady state and the solution is null velocity and flat water surface.\\

In SWASHES we consider the following initial conditions:
\[h+z=\max(0.1,z)\;\text{m and}\;q=0\;\text{m}^2\text{/s}\]
and the boundary conditions
\begin{equation*}
\left\{\begin{array}{l}
h=0.1\;\text{m},\\
q=0\;\text{m}^2\text{/s}.\end{array}\right.
 \end{equation*}


\subsubsection{Subcritical flow}
\label{Bump-Subcritical}
After testing the two steady states at rest, the user can increase the difficulty with dynamical steady states.
In the case of a subcritical flow, using \eqref{eq:Bernoulli}, the water height is given by the resolution of
 \begin{equation*}
h(x)^3+\left(z(x)-\dfrac{{q_0}^2}{2gh_{L}^2}-h_{L}\right)h(x)^2+\dfrac{q_0^2}{2g}=0,\;\forall x\in[0,L],
\end{equation*}
 where $h_{L} = h(x=L)$.\\

Thanks to the Bernoulli relation \eqref{eq:Bernoulli}, we can notice that the water height is constant when the topography is constant, decreases (respectively increases) when the bed slope 
 increases (resp. decreases). The water height reaches its minimum at the top of the bump \cite{Goutal97}.

In SWASHES, the initial conditions are
\[h+z=2\;\text{m and}\;q=0\;\text{m}^2\text{/s}\]
and the boundary conditions are chosen as
\begin{equation*}
  \left\{\begin{array}{l}
    \text{upstream:}\;q=4.42\;\text{m}^2\text{/s},\\
    \text{downstream:}\;h=2\;\text{m}.
  \end{array}\right.
\end{equation*}


\subsubsection{Transcritical flow without shock}
\label{Bump-Transcritical-NoShock}
In this part, we consider the case of a transcritical flow, without shock \cite{Goutal97}.
Again thanks to \eqref{eq:Bernoulli}, we can express the water height as the solution of
 \begin{equation*}
h(x)^3+\left(z(x)-\dfrac{{q_0}^2}{2gh_{c}^2}-h_{c}-z_{M}\right)h(x)^2+\dfrac{q_0^2}{2g}=0,\quad \forall x\in[0,L],
\end{equation*}
 where $z_{M}=\max_{x\in[0,L]}{z}$ and $h_c$ is the corresponding water height.\\

The flow is fluvial upstream and becomes torrential at the top of the bump.

Initial conditions can be taken equal to
\[h+z=0.66\;\text{m and}\;q=0\;\text{m}^2\text{/s}\]
and for the boundary conditions
\begin{equation*}
  \left\{\begin{array}{l}
    \text{upstream:}\;q=1.53\;\text{m}^2\text{/s},\\
    \text{downstream:}\;h=0.66\;\text{m while the flow is subcritical}.
  \end{array}\right.
\end{equation*}


\subsubsection{Transcritical flow with shock}
\label{Bump-Transcritical-WithShock-NoFriction}
If there is a shock in the solution \cite{Goutal97}, using \eqref{eq:Bernoulli}, the water height is given by the resolution of
 \begin{equation}
\left\{\begin{array}{l  l}
h(x)^3+\left(z(x)-\dfrac{{q_0}^2}{2gh_{c}^2}-h_{c}-z_{M}\right)h(x)^2+\dfrac{q_0^2}{2g}=0 & \mbox{ for } x < x_{shock}, \\[0.3cm]
h(x)^3+\left(z(x)-\dfrac{{q_0}^2}{2gh_{L}^2}-h_{L}\right)h(x)^2+\dfrac{q_0^2}{2g}=0 & \mbox{ for } x > x_{shock},\\[0.3cm]
 {q_0}^2\left(\dfrac{1}{{h_1}}-\dfrac{1}{{h_2}}\right)+\dfrac{g}{2}\left({h_1}^2-{h_2}^2\right)=0.&
 \end{array} \right.\label{eq:bump-shock}
\end{equation}
In these equalities, $z_{M}=\max_{x\in[0,L]}{z}$, $h_c$ is the corresponding water height, $h_{L} = h(x=L)$ 
and $h_{1} = h(x_{shock}^-)$, $h_{2}= h(x_{shock}^+)$ are the water heights upstream
and downstream respectively.
The shock is located thanks to the third relation in system \eqref{eq:bump-shock}, which is a Rankine-Hugoniot's relation.\\

As for the previous case, the flow becomes supercritical at the top of the bump but it becomes again fluvial after a hydraulic jump.

One can choose for initial conditions
\[h+z=0.33\;\text{m and}\;q=0\;\text{m}^2\text{/s}\]
and the following boundary conditions
\begin{equation*}
  \left\{\begin{array}{l}
    \text{upstream:}\;q=0.18\;\text{m}^2\text{/s},\\
    \text{downstream:}\;h=0.33\;\text{m}.
  \end{array}\right.
\end{equation*}

We can find a generalisation of this case with a friction term in Hervouet's work \cite{Hervouet07}.


\subsection{Mac Donald's type 1D solutions}
\label{MacDo-1D}

Following the lines of \cite{MacDonald96,MacDonald97}, we give here some steady state solutions of system~\eqref{SaintVenant1D} with
 varying topography and friction term (from \cite{Vo08, Delestre10b}). Rain and diffusion are not considered
 ($R=0~\text{m}/\text{s}$, $\mu=0~\text{m}^2/\text{s}$), so the steady states system~\eqref{eq:general-steady-SV} reduces to
\begin{equation}
 \partial_x z=\left(\dfrac{q^2}{g{h}^3}-1\right)\partial_x {h}-S_f.\label{eq:steady-SV}
\end{equation}
From this relation, one can make as many solutions as required. In this section, we present some of them that are obtained
for specific values of the length $L$ of the domain, and for fixed parameters (such as the friction law and its coefficient).
The water height profile and the discharge are given, and we compute the corresponding topographies
solving equation~\eqref{eq:steady-SV}.
We have to mention that there exists another approach, classical in hydraulics. It consists in considering a given topography and a discharge.
From these, the steady-state water height is deduced thanks to equation~\eqref{eq:steady-SV} and to the classification of water surface profiles
(see among others \cite{Chow59} and \cite{Henderson66}). Solutions obtained using this approach may be found in \cite{Zhou99} for example.
Finally, note that some simple choices for the free surface may lead to exact analytic solutions.
 
The solutions given in this section are more intricate than the ones of the previous section, as the topography can vary
near the boundary. Consequently they give a better validation of the boundary conditions.
If $S_{f} \not = 0$ (we have friction at the bottom), the following solutions can prove if the
friction terms are coded in order to satisfy the steady states. 
\begin{remark}
\label{rem:raff}
All these solutions are given by the numerical resolution of an equation.
So, the space step should be small enough to have a sufficiently precise solution.
It means that the space step used to get these solutions should be smaller than the 
space step of the code to be validated.
\end{remark}


\subsubsection{Long channel: 1000~m}
\label{McDo-1D-1000m}


\paragraph{Subcritical case}
\label{McDo-1D-1000m-Subcritical}

\noindent We consider a 1000~m long channel with a discharge of $q=2\;\text{m}^2/\text{s}$ \cite{Vo08}. The flow is
 constant at inflow and the water height is prescribed at outflow, with the following values:
\begin{equation*}
 \left\{\begin{array}{l}
   \text{upstream:}\;q=2\;\text{m}^2/\text{s},\\
   \text{downstream:}\;h=h_{ex}(1000).
 \end{array}\right.
\end{equation*}
The channel is initially dry, \emph{i.e.} initial conditions are
\begin{equation*}
h=0 \;\text{m}\;\text{ and }\;q=0\;\text{m}^2/\text{s}.
\end{equation*}
The water height is given by
\begin{equation}
h_{ex}(x)=\left(\dfrac{4}{g}\right)^{1/3}\left(1+\dfrac{1}{2}\exp\left(-16\left(\dfrac{x}{1000}-\dfrac{1}{2}\right)^2\right)\right).\label{eqfluvfluv}
\end{equation}
We remind the reader that $q=2\;\text{m}^2/\text{s}$ on the domain and that the topography is calculated iteratively
thanks to \eqref{eq:steady-SV}.
We can consider the two friction laws explained in the introduction, with the coefficients
$n = 0.033$~m\textsuperscript{-1/3}s for Manning's and $f=0.093$ for Darcy-Weisbach's.\\
Under such conditions, we get a subcritical steady flow.


\paragraph{Supercritical case}
\label{McDo-1D-1000m-Supercritical}
We still consider a 1000~m long channel, but with a constant discharge
 $q=2.5 \;\text{m}^2/\text{s}$ on the whole domain  \cite{Delestre10b}. The flow is supercritical both at inflow and at outflow, 
thus we consider the following boundary conditions:
\begin{equation*}
 \left\{\begin{array}{l}
   \text{upstream:}\;q=2.5\;\text{m}^2/\text{s} \;\text{ and }\; h=h_{ex}(0),\\
   \text{downstream: free.}
 \end{array}\right.
\end{equation*}
 The initial conditions are a dry channel
\begin{equation*}
h=0 \;\text{m}\;\text{ and }\;q=0\;\text{m}^2/\text{s}.
\end{equation*}
With the water height given by
\begin{equation}
h_{ex}(x)=\left(\dfrac{4}{g}\right)^{1/3}\left(1-\dfrac{1}{5}\exp\left(-36\left(\dfrac{x}{1000}-\dfrac{1}{2}\right)^2\right)\right)\label{eqtortor}
\end{equation}
and the friction coefficients equal to $n = 0.04$~m\textsuperscript{-1/3}s for Manning's and $f= 0.065 $ for Darcy-Weisbach's friction law,
the flow is supercritical.


\paragraph{Subcritical-to-supercritical case}
\label{McDo-1D-1000m-SubSupercritical}

The channel is 1000~m long and the discharge at equilibrium is $q=2 \;\text{m}^2/\text{s}$ \cite{Vo08}. The flow is subcritical upstream
 and supercritical downstream, thus we consider
the following boundary conditions:
\begin{equation*}
 \left\{\begin{array}{l}
   \text{upstream:}\;q=2\;\text{m}^2/\text{s},\\
   \text{downstream: free.}
  \end{array}\right.
\end{equation*}
As initial conditions, we consider a dry channel
\begin{equation*}
h=0 \;\text{m}\;\text{ and }\;q=0\;\text{m}^2/\text{s}.
\end{equation*}
In this configuration, the water height is
\begin{equation*}
h_{ex}(x)=\left\{\begin{array}{ll}
        \left(\dfrac{4}{g}\right)^{1/3}\left(1-\dfrac{1}{3}\tanh\left(3\left(\dfrac{x}{1000}-\dfrac{1}{2}\right)\right)\right) & \mbox{ for } 0\text{ m}\leq x \leq 500\text{ m}, \\
        \left(\dfrac{4}{g}\right)^{1/3}\left(1-\dfrac{1}{6}\tanh\left(6\left(\dfrac{x}{1000}-\dfrac{1}{2}\right)\right)\right) & \mbox{ for } 500\text{ m}< x \leq 1000\text{ m},
            \end{array}\right.
\end{equation*}
with a friction coefficient $n = 0.0218$~m\textsuperscript{-1/3}s (resp. $f= 0.042 $) for the Manning's (resp. the Darcy-Weisbach's) law. Thus we
 get a transcritical flow (from fluvial to torrential via a transonic point).


\paragraph{Supercritical-to-subcritical case}
\label{McDo-1D-1000m-SuperSubcritical}

As in the previous cases, the domain is 1000~m long and the discharge is $q=2 \;\text{m}^2/\text{s}$ \cite{Vo08}. The boundary conditions are
 a torrential inflow and a fluvial outflow:
\begin{equation*}
 \left\{\begin{array}{l}
   \text{upstream:}\;q=2\;\text{m}^2/\text{s} \text{ and } h=h_{ex}(0),\\
   \text{downstream:}\;h=h_{ex}(1000).
 \end{array}\right.
\end{equation*}
At time $t=0$ s, the channel is initially dry
\begin{equation*}
h=0 \;\text{m}\;\text{ and }\;q=0\;\text{m}^2/\text{s}.
\end{equation*}
The water height is defined by the following discontinuous function
\begin{equation*}
h_{ex}(x)=\left\{\begin{array}{ll}
              \left(\dfrac{4}{g}\right)^{1/3}\left(\dfrac{9}{10}-\dfrac{1}{6}\exp\left(-\dfrac{x}{250}\right)\right) 
             & \hspace{-5cm}  \mbox{ for } 0\text{ m}\leq x \leq 500\text{ m},  \\[0.7cm]
      \left(\dfrac{4}{g}\right)^{1/3}\left(1+\displaystyle{\sum\limits_{ k=1}^3} a_k\exp\left(-20k\left(\dfrac{x}{1000} -\dfrac{1}{2}\right)\right)
          +\dfrac{4}{5}\exp\left(\dfrac{x}{1000}-1\right)\right)&\\
             & \hspace{-5cm}  \mbox{ for } 500\text{ m}\leq x \leq 1000\text{ m},
             \end{array}\right.
\end{equation*}
with  \(a_1=-0.348427\), \(a_2=0.552264\), \(a_3=-0.55558\).\\
The friction coefficients are $n = 0.0218$~m\textsuperscript{-1/3}s for the Manning's law and $f= 0.0425 $ for the Darcy-Weisbach's law. The steady
 state solution is supercritical upstream and becomes subcritical through a hydraulic jump located at $x=500$~m.


\subsubsection{Short channel: 100~m}
\label{McDo-1D-100m}

In this part, the friction law we consider is the Manning's law.
Generalization to other classical friction laws is straightforward.


\paragraph{Case with smooth transition and shock}
\label{McDo-1D-100m-Subcritical}

The length of the channel is $100$ m and the discharge at steady states is $q=2 \;\text{m}^2/\text{s}$ \cite{Vo08}. The flow
 is fluvial both upstream and downstream, the boundary conditions are fixed as follows
\begin{equation*}
 \left\{\begin{array}{l}
   \text{upstream:}\;q=2\;\text{m}^2/\text{s},\\
   \text{downstream:}\;h=h_{ex}(100).
 \end{array}\right.
\end{equation*}
To have a case including two kinds of flow (subcritical and supercritical)
and two kinds of transition (transonic and shock), we consider a channel filled with water, \emph{i.e.}
\begin{equation*}
h(x)=\max(h_{ex}(100)+z(100)-z(x),0)\text{ and }\;q=0\;\text{m}^2/\text{s}.
\end{equation*}
The water height function has the following formula
\begin{equation*}
h_{ex}(x)=\left\{\begin{array}{ll}
              \left(\dfrac{4}{g}\right)^{1/3}\left(\dfrac{4}{3}-\dfrac{x}{100}\right)-\dfrac{9x}{1000}\left(\dfrac{x}{100}-\dfrac{2}{3}\right)&\\
              &\hspace{-4cm}\mbox{ for } 0\text{ m}\leq x \leq \dfrac{200}{3}\approx 66.67\text{ m},\\[0.7cm]
    \left(\dfrac{4}{g}\right)^{1/3}\left(a_1 \left(\dfrac{x}{100}-\dfrac{2}{3}\right)^4
+a_1 \left(\dfrac{x}{100}-\dfrac{2}{3}\right)^3\right. & \\
\left. \qquad -a_2 \left(\dfrac{x}{100}-\dfrac{2}{3}\right)^2+a_3\left(\dfrac{x}{100}-\dfrac{2}{3}\right)+a_4\right)&\\
              &\hspace{-4cm} \mbox{ for } \dfrac{200}{3}\approx 66.67\text{ m}\leq x \leq 100\text{ m},
             \end{array}\right.
\end{equation*}
with $a_1=0.674202$, $a_2=21.7112$, $a_3=14.492$ et $a_4=1.4305$.\\
In this case, the Manning's friction coefficient is $n = 0.0328$~m\textsuperscript{-1/3}s, the inflow is subcritical,
becomes supercritical via a sonic point, and, through a shock (located at $x=200/3 \approx 66.67$~m), becomes subcritical again.


\paragraph{Supercritical case}
\label{McDo-1D-100m-Supercritical}

The channel we consider is still 100~m long and the equilibrium discharge is $q=2 \;\text{m}^2/\text{s}$ \cite{Delestre10b}. The flow is torrential
at the bounds of the channel, thus the boundary conditions are
\begin{equation*}
 \left\{\begin{array}{l}
   \text{upstream:}\;q=2\;\text{m}^2/\text{s} \text{ and } h=h_{ex}(0),  \\
   \text{downstream: free.}
 \end{array}\right.
\end{equation*}
As initial conditions, we consider an empty channel which writes
\begin{equation*}
h=0 \;\text{m}\;\text{ and }\;q=0\;\text{m}^2/\text{s}.
\end{equation*}
The water height is given by
\begin{equation*}
 h_{ex}(x)=\left(\dfrac{4}{g}\right)^{1/3}\left(1-\dfrac{1}{4}\exp\left(-4\left(\dfrac{x}{100}-\dfrac{1}{2}\right)^2\right)\right)
\end{equation*}
and the friction coefficient is $n = 0.03$~m\textsuperscript{-1/3}s (for the Manning's law). The flow is entirely torrential.


\paragraph{Subcritical-to-supercritical case}
\label{McDo-1D-100m-SubSupercritical}

A $100$ m long channel has a discharge of $q=2 \;\text{m}^2/\text{s}$ \cite{Vo08}.
The flow is fluvial at inflow and torrential at outflow with following boundary conditions
\begin{equation*}
 \left\{\begin{array}{l}
   \text{upstream:}\;q=2\;\text{m}^2/\text{s},\\
   \text{downstream: free.}
 \end{array}\right.
\end{equation*}
As in the subcritical case, the initial condition is an empty channel with a puddle downstream
\begin{equation*}
h(x)=\max(h_{ex}(100)+z(100)-z(x),0)\text{ and }\;q=0\;\text{m}^2/\text{s},
\end{equation*}
and the water height is
\begin{equation*}
h_{ex}(x)=\left(\dfrac{4}{g}\right)^{1/3}\left(1-\dfrac{\left(x-50\right)}{200}+\dfrac{\left(x-50\right)^2}{30000}\right).
\end{equation*}
The Manning's friction coefficient for the channel is $n = 0.0328$~m\textsuperscript{-1/3}s. 
We get a transcritical flow: subcritical upstream and
 supercritical downstream.


\subsubsection{Long channel: 5000~m, periodic and subcritical}
\label{McDo-1D-5000m-Periodic-Subcritical}

For this case, the channel is much longer than for the previous ones: 5000~m, but the discharge at equilibrium is still
 $q=2 \;\text{m}^2/\text{s}$ \cite{Vo08}. Inflow and outflow are both subcritical. The boundary conditions are taken as:
\begin{equation*}
 \left\{\begin{array}{l}
   \text{upstream:}\;q=2\;\text{m}^2/\text{s},\\
   \text{downstream:}\;h=h_{ex}(5000).
 \end{array}\right.
\end{equation*}
We consider a dry channel with a little lake at rest downstream as initial conditions:
\begin{equation*}
h(x)=\max(h_{ex}(5000)+z(5000)-z(x),0) \text{ and }q=0\;\text{m}^2/\text{s}.
\end{equation*}
We take the water height at equilibrium as a periodic function in space, namely
\begin{equation*}
h_{ex}(x)=\dfrac{9}{8}+\dfrac{1}{4}\sin\left(\dfrac{\pi x}{500}\right)
\end{equation*}
and the Manning's constant is $n = 0.03$~m\textsuperscript{-1/3}s. 
We get a subcritical flow. As the water height is periodic, the associated topography
(solution of Equation \eqref{eq:steady-SV}) is periodic as well:
we get a periodic configuration closed to the ridges-and-furrows configuration. Thus this case is interesting for
 the validation of numerical methods for overland flow simulations on agricultural fields.


\subsection{Mac Donald's type 1D solutions with rain}
\label{McDo-1Drain}

In this section, we consider the Shallow Water system~\eqref{SaintVenant1D} with rain (but without viscosity: $\mu=0~\text{m}^2/\text{s}$)
 at steady states \cite{Vo08}. 
The rain intensity is constant, equal to $R_{0}$. The rain is uniform on the domain $[0,L]$.
Under these conditions, if we denote by $q_{0}$ the discharge value at inflow $q(t,0)=q_0$, we have:
\begin{equation}
 q(x)=q_0+xR_0, \text{ for } 0\leq x \leq L.\label{macdoDebitP}
\end{equation}
The solutions are the same as in section~\ref{MacDo-1D}, except that the discharge is given by~\eqref{macdoDebitP}. 
But the rain term modifies the expression of the topography through a new rain term as written in~\eqref{eq:general-steady-SV}. More precisely, 
Equation~\eqref{eq:steady-SV} is replaced by
\begin{equation*}
 \partial_x z=\left(\dfrac{q^2}{gh^3}-1\right)\partial_x h-\dfrac{2qR_0}{gh^2}-S_f. \label{bressePluie}
\end{equation*}

These solutions allow the validation of the numerical treatment of the rain.
Remark~\ref{rem:raff}, mentioned in the previous section, applies to these solutions too.


\subsubsection{Subcritical case for a long channel}
\label{McDo-1Drain-Subcritical}
For a 1000~m long channel, we consider a flow which is fluvial on the whole domain. Thus we impose the following boundary
 conditions:
\begin{equation*}
 \left\{\begin{array}{l}
   \text{upstream:}\;q=q_0,\\
   \text{downstream:}\;h=h_{ex}(1000),
 \end{array}\right.
\end{equation*}
with the initial conditions 
$$ h = 0~\text{m} \mbox{ and } q=0~\text{m}^2/\text{s},$$
where $h_{ex}$ is the water height at steady state given by \eqref{eqfluvfluv}.

In SWASHES, for the friction term, we can choose either Manning's law with $n = 0.033$~m\textsuperscript{-1/3}s or Darcy-Weisbach's law with $f=0.093$, 
 the discharge $q_{0}$ is fixed at $1~\text{m}^2/\text{s}$ and the rain intensity is $R_{0} = 0.001$~m/s.


\subsubsection{Supercritical case for a long channel}
\label{McDo-1Drain-Supercritical}
The channel length remains unchanged (1000~m), but, as the flow is supercritical, the boundary conditions are
\begin{equation*}
 \left\{\begin{array}{l}
   \text{upstream:}\;q = q_{0} \mbox{ and } h=h_{ex}(0),\\
   \text{downstream: free.}
 \end{array}\right.
\end{equation*}
At initial time, the channel is dry
$$ h = 0~\text{m} \mbox{ and } q=0~\text{m}^2/\text{s}.$$
At steady state, the formula for the water height is \eqref{eqtortor}. 

From a numerical point of view, for this case, we recommend that the rain does not start at the initial time.
 The general form of the recommended rainfall event is
\begin{equation*}
R(t)=\left\{\begin{array}{ll}
         0 \text{ m/s}    &  \text{if}\;t<t_{R},\\
  R_{0}   &  \text{else},
        \end{array}\right.
\end{equation*}
with $t_{R} = 1500$~s. Indeed, this allows to get two successive steady states: for the first one the discharge is constant
 in space $q_0$ and for the second one the discharge is \eqref{macdoDebitP} with the chosen height profile \eqref{eqtortor}.

In SWASHES, we have a friction coefficient $n = 0.04$~m\textsuperscript{-1/3}s for Manning's law, $f = 0.065$ for Darcy-Weisbach's law. Inflow
 discharge is $q_0=2.5\;\text{m}^2/\text{s}$ and $R_{0} = 0.001$~m/s.


\subsection{Mac Donald's type 1D solutions with diffusion}
\label{McDo-1Ddiffusion}

Following the lines of \cite{MacDonald96,MacDonald97}, Delestre and Marche, in \cite{Delestre10}, proposed new analytic solutions with a diffusion
 term (with $R=0$~m/s). To our knowledge, these are the only analytic solutions available in the litterature with a
 diffusion source term.
In \cite{Delestre10}, the authors considered system \eqref{SaintVenant1D} with the source terms derived in \cite{Marche07},
 \emph{i.e.}:
\begin{equation*}
 S_f=\dfrac{1}{g}\left(\dfrac{\alpha_0(h)}{h}u+\alpha_1(h)|u|u\right)\quad\text{and}\quad \mu=4\mu_h
\end{equation*}
 with
\begin{equation*}
 \alpha_0(h)=\dfrac{k_l}{1+\dfrac{k_l h}{3 \mu_v}} \quad\text{and}\quad
 \alpha_1(h)=\dfrac{k_t}{\left(1+\dfrac{k_l h}{3 \mu_v}\right)^2}
\end{equation*}
where $\mu_v$ [T] (respectively $\mu_h$ [L$^2$/T]) is the vertical (resp. the horizontal) eddy viscosity and $k_l$ [T/L] (resp. $k_t$ [1/L]) the laminar
 (resp. the turbulent) friction coefficient.
At steady states we recover \eqref{eq:general-steady-SV} with $R=0~\text{m}/\text{s}$, or again
\begin{equation}
 \left\{\begin{array}{l}
         q=q_0,\\
  \partial_x z=\dfrac{1}{gh}\left(\dfrac{q^2}{h^2}-gh\right)\partial_x h
-S_f(h,q)+\dfrac{\mu}{gh^2}\left(-q\partial_{xx}h+\dfrac{q}{h}(\partial_x h)^2\right).
        \end{array}\right.\label{Bresse-diff}
\end{equation}
 As for the previous cases, the topography is evaluated thanks to the momentum equation of \eqref{Bresse-diff}.

 These cases allow for the validation of the diffusion source term treatment. These solutions may be easily adapted to
 Manning's and Darcy-Weisbach's friction terms. Remark~\ref{rem:raff} applies to these solutions too.

 In \cite{Delestre10}, the effect of $\mu_h$, $\mu_v$, $k_t$ and $k_l$ is studied by using several values. In what follows,
 we present only two of these solutions: a subcritical flow and a supercritical flow.


\subsubsection{Subcritical case for a long channel}
\label{McDo-1Ddiffusion-Subcritical}
A $1000$~m long channel has a discharge of $q=1.5\;\text{m}^2/\text{s}$. The flow is fluvial at both channel boundaries,
 thus the boundary conditions are:
\begin{equation*}
 \left\{\begin{array}{l}
  \text{upstream:}\;q=1.5\;\text{m}^2/\text{s},\\
  \text{downstream:}\;h=h_{ex}(1000).
 \end{array}\right.
\end{equation*}
The channel is initially dry
\begin{equation*}
 h=0\;\text{m}\quad\text{and}\quad q=0\;\text{m}^2/\text{s}
\end{equation*}
and the water height at steady state ($h_{ex}$) is the same as in section~\ref{McDo-1D-1000m-Subcritical}.

In SWASHES, the parameters are: $k_t = 0.01$, $k_l = 0.001$, $\mu_v = 0.01$ and $\mu_h = 0.001$.


\subsubsection{Supercritical case for a long channel}
\label{McDo-1Ddiffusion-Supercritical}
We still consider a $1000$~m long channel, with a constant discharge
 $q=2.5\;\text{m}^2/\text{s}$ on the whole domain. Inflow and outflow are both torrential, thus we choose the following
 boundary conditions:
\begin{equation*}
 \left\{\begin{array}{l}
    \text{upstream:}\;q=2.5\;\text{m}^2/\text{s}\;\text{and}\; h=h_{ex}(0),\\
    \text{downstream: free.}
 \end{array}\right.
\end{equation*}
We consider a dry channel as initial condition:
\begin{equation*}
 h=0\;\text{m}\;\text{and}\;q=0\;\text{m}^2/\text{s}.
\end{equation*}
The water height $h_{ex}$ at steady state is given by function \eqref{eqtortor}.

In SWASHES, we have: $k_t = 0.005$, $k_l = 0.001$, $\mu_v = 0.01$ and $\mu_h = 0.1$.


\subsection{Mac Donald pseudo-2D solutions}
\label{McDo-pseudo2D}

In this section, we give several analytic solutions for the pseudo-2D Shallow-Water system. 
This system can be considered as an intermediate between the one-dimensional and the two-dimensional models. 
More precisely, these equations model a flow in a rectilinear three-dimensional channel
with the quantities averaged not only on the vertical direction but also on the width of the channel.
For the derivation, see for example \cite{Goutal10}.
Remark~\ref{rem:raff}, mentioned for Mac Donald's type 1D solutions, applies to these pseudo-2D solutions too.

We consider six cases for non-prismatic channels introduced in \cite{MacDonald96}.
These channels have a variable slope and their width is also variable in space.
More precisely, each channel is determined through the definition of the bottom width $B$ (as a function of the space variable $x$)
and the slope of the boundary $Z$ (Figure~\ref{figMacDo2d}). The bed slope is an explicit function of the water height,
detailed in the following. 

\begin{figure}[htbp]
\begin{center}
\includegraphics[width = 0.9 \textwidth]{MacDonald_section} 
\caption{Notations for the Mac Donald pseudo-2D solutions}
\label{figMacDo2d}
\end{center}
\end{figure}

\begin{table}[htbp]
\begin{center}
  \begin{tabular}{|p{0.35\textwidth}|c|c|c|c|c|}
    \hline
    Solution & $B$ (m) & Z (m)  & $L$ (m) &  $h_{in}$ (m) & $h_{out}$ (m) \\[0.1cm]
    \hline
    \hline
     Subcritical flow in a short domain& $B_{1}(x)$ & 0 & 200  & & 0.902921 \\
    \hline
    Supercritical flow  in a short domain & $B_{1}(x)$ & 0 & 200 & 0.503369 & \\
    \hline
    Smooth transition  in a short domain& $B_{1}(x)$ & 0 & 200 & & \\
    \hline
    Hydraulic jump  in a short domain& $B_{1}(x)$ & 0 & 200 &  0.7 & 1.215485 \\
    \hline
    Subcritical flow in a long domain & $B_{2}(x)$ & 2 & 400 &  & 0.904094 \\
    \hline
    Smooth transition followed by a hydraulic jump in a long domain& $B_{2}(x)$ & 2 & 400 &  & 1.2 \\
    \hline
  \end{tabular}\\[1.2\baselineskip]
\end{center}
\caption{Main features of the cases of pseudo-2D channels}
\label{tab:pseudo2D}
\end{table}

\begin{figure}
\begin{center}
\includegraphics[width = 0.49 \textwidth]{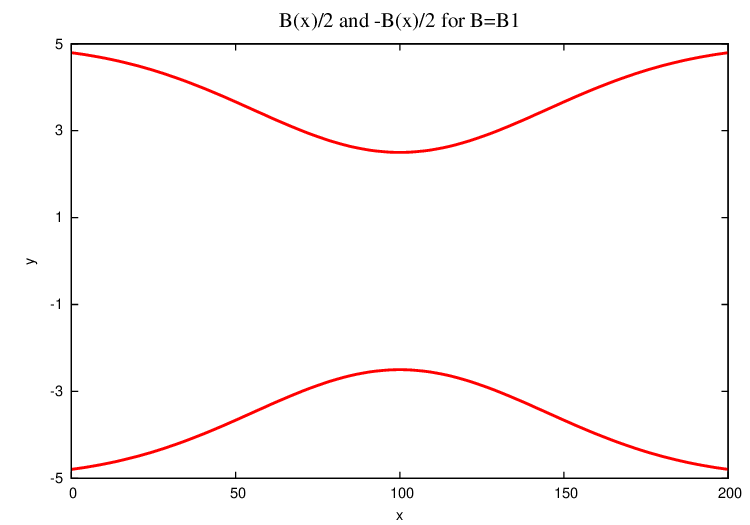} \includegraphics[width = 0.49 \textwidth]{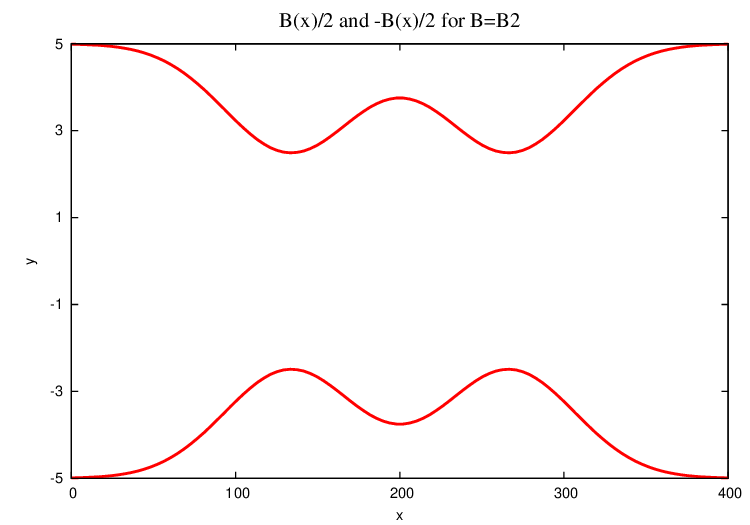} 
\caption{Functions $B_{1}$ and $B_{2}$ for the shape of the channel}
\label{figB}
\end{center}
\end{figure}

\noindent
The features of these cases are summarized in table~\ref{tab:pseudo2D}.
In this table, the functions for the bed shape are:
\begin{equation*}
\begin{array}{c}
B_{1}(x)=10 - 5 \exp \left( -10 \left( \dfrac{x}{200} - \dfrac{1}{2} \right)^{2}\right)  \  \mbox{for } 0\mbox{ m } \leq x \leq L = 200 \mbox{ m }, \\
B_{2}(x)=10 - 5 \exp \left( -50 \left( \dfrac{x}{400} - \dfrac{1}{3} \right)^{2}\right) - 5 \exp \left( -50 \left( \dfrac{x}{400} - \dfrac{2}{3} \right)^{2}\right)
  \\  \mbox{for } 0\mbox{ m } \leq x \leq L = 400 \mbox{ m },
\end{array}
\end{equation*}
(Figure~\ref{figB}) and $h_{in}$ (resp. $h_{out}$) is the water height at the inflow (resp. outflow).  \\

\noindent
In each case, the Manning's friction coefficient is $n = 0.03$~m\textsuperscript{-1/3}s, the discharge $q$ is taken equal to
$20\text{ m}^3 \text{s}^{-1}$, the slope of the topography is given by:
\begin{multline*}
S_{0}(x)=\left(1-\dfrac{q^{2}(B(x)+2 Z h(x))}{g 
                        \,h(x)^{3} (B(x)+Z h(x))^{3}}\right) h'(x)
  + q^{2} n^{2} \dfrac{\left(B(x)+2 h(x) \sqrt{1+Z^{2}}\right)^{{4}/{3}}}{h(x)^{{10}/{3}}\left(B(x)+Z h(x)\right)^{{10}/{3}}}\\
  - \dfrac{q^{2}B'(x)}{g 
            \,h(x)^{2} (B(x) + Z h(x))^{3}}
\end{multline*}
where $h$ is the water height and the topography is defined as $z(x)= \displaystyle \int_x^L S_{0}(X)\;dX$.

\begin{remark}
When programming these formulae, we noted a few typos in \cite{MacDonald96}, in the expression of $S_{0}(x)$,
of $\phi$ in the hydraulic jump case (\textsection~\ref{McDo-pseudo2D-Jump}) and in the formulation of $h$ in $[0;120]$
in the solution for the smooth transition followed by a hydraulic jump (\textsection~\ref{McDo-pseudo2D-SmoothJump}).
\end{remark}

\begin{remark}
We recall that the following analytic solutions are solutions of the pseudo-2D Shallow-Water system. 
This is the reason why, in this section, $h_{ex}$ does not depend on $y$. 
\end{remark}

\subsubsection{Subcritical flow in a short domain}
\label{McDo-pseudo2D-Subcritical}
In the case of a subcritical flow in a short domain, as in the three that follow, the cross section of the channel is rectangular, the bottom is given by the function $B_{1}$
and the length $L=200$~m.

In this current case, the flow is fixed at inflow and the water height is prescribed at outflow. We have the following boundary conditions:
\begin{equation*}
 \left\{\begin{array}{l}
    \text{upstream:}\;q=20\text{ m}^3 \text{s}^{-1},\\
    \text{downstream:}\;h=h_{out}.
 \end{array}\right.
\end{equation*}

The channel is initially dry, with a little puddle downstream, \emph{i.e.} initial conditions are:
\begin{equation*}
 h(x,y)= \max(0, h_{out} + z(200,y) - z(x,y))\;\text{and}\;q=0\;\text{m}^3/\text{s}.
\end{equation*}

If we take the mean water height 
$$h_{ex}(x)=0.9 + 0.3 \exp\left(-20 \left( \frac{x}{200} - \frac{1}{2} \right)^{2} \right), $$
the flow stays subcritical in the whole domain of length $L = 200$~m.

\subsubsection{Supercritical flow in a short domain}
\label{McDo-pseudo2D-Supercritical}

In this case, the flow and the water height are fixed at inflow. We have the following boundary conditions:
\begin{equation*}
  \left\{\begin{array}{l}
     \text{upstream:}\;q=20\text{ m}^3 \text{s}^{-1}\; \text{and}\; h=h_{in},\\
     \text{downstream: free.}
  \end{array}\right.
\end{equation*}

The channel is initially dry, \emph{i.e.} initial conditions are:
\begin{equation*}
  h=0\text{ m},\;\text{and}\;q=0\;\text{m}^3/\text{s}.
\end{equation*}

If we consider the mean water height
$$h_{ex}(x)=0.5 + 0.5 \exp\left(-20 \left( \frac{x}{200} - \frac{1}{2} \right)^{2} \right),$$
in a channel of length $L = 200$~m with the $B_{1}$ shape and vertical boundary, the flow is supercritical.

\subsubsection{Smooth transition in a short domain}
\label{McDo-pseudo2D-Smooth}

In the case of smooth transition in a short domain, the flow is fixed at the inflow. We have the following boundary conditions:
\begin{equation*}
   \left\{\begin{array}{l}
      \text{upstream:}\;q=20\text{ m}^3 \text{s}^{-1},\\
      \text{downstream: free.}
   \end{array}\right.
\end{equation*}

The channel is initially dry, \emph{i.e.} initial conditions are:
\begin{equation*}
  h=0\text{ m},\;\text{and}\;q=0\;\text{m}^3/\text{s}.
\end{equation*}

The channel is the same as in the previous cases, with a mean water height given by
$$h_{ex}(x)=1 - 0.3 \tanh \left(4 \left( \frac{x}{200} - \frac{1}{3} \right) \right).$$
Under these conditions, the flow is first subcritical and becomes supercritical.

\subsubsection{Hydraulic jump  in a short domain}
\label{McDo-pseudo2D-Jump}

In the case of a hydraulic jump  in a short domain, the flow discharge is fixed at the inflow and the water height is prescribed at both inflow and outflow. We have the following boundary conditions:
\begin{equation*}
   \left\{\begin{array}{l}
      \text{upstream:}\;q=20\text{ m}^3 \text{s}^{-1}\; \text{and}\; h=h_{in},\\
      \text{downstream:}\;h=h_{out}.
   \end{array}\right.
\end{equation*}

The channel is initially dry, with a little puddle downstream, \emph{i.e.} initial conditions are:
\begin{equation*}
 h(x,y)= \max(0, h_{out} + z(200,y) - z(x,y))\;\text{and}\;q=0\;\text{m}^3/\text{s}.
\end{equation*}

We choose the following expression for the mean water height:
$$h_{ex}(x)=0.7 + 0.3 \left(\exp\left(\frac{x}{200}\right) - 1 \right)\quad \mbox{ for } 0\mbox{ m} \leq  x  \leq 120 \mbox{ m},$$
and
$$h_{ex}(x)=\exp\left(-p (x - x^{\star}) \right) \sum_{i=0}^{M} k_{i} \left( \frac{x-x^{\star}}{x^{\star\star}-x^{\star}} \right)^{i} + \phi(x)
  \quad \mbox{ for } 120\mbox{ m} \leq  x  \leq 200 \mbox{ m},$$
with:\\[0.1cm]
\begin{tabular}{lll}
$\bullet$ $x^{\star}=120$~m,&
$\bullet$ $x^{\star\star}=200$~m,&
$\bullet$ $M=2$,\\
$\bullet$ $k_{0}=-0.154375$,&
$\bullet$ $k_{1}=-0.108189$,&
$\bullet$ $k_{2}=-2.014310$,\\
$\bullet$ $p=0.1$,&
$\bullet$ $\phi(x)=1.5 \exp \left(0.1 \left( \dfrac{x}{200} - 1 \right) \right)$.
\end{tabular}\\[0.1cm]
We obtain a supercritical flow that turns into a subcritical flow through a hydraulic jump.

\subsubsection{Subcritical flow in a long domain}
\label{McDo-pseudo2D-Subcritical2}

From now on, the length of the domain is $L = 400$~m, the boundaries of the channel are given by $B_{2}$ and the cross sections
are isoscele trapezoids.

In this case, the flow is fixed at the inflow and the water height is prescribed at the outflow. We have the following boundary conditions:
\begin{equation*}
  \left\{\begin{array}{l}
     \text{upstream:}\;q=20\text{ m}^3 \text{s}^{-1},\\
     \text{downstream:}\;h=h_{out}.
  \end{array}\right.
 \end{equation*}

The channel is initially dry, with a little puddle downstream, \emph{i.e.} initial conditions are:
\begin{equation*}
 h(x,y)= \max(0, h_{out} + z(400,y) - z(x,y))\;\text{and}\;q=0\;\text{m}^3/\text{s}.
\end{equation*}

Considering the mean water height
$$h_{ex}(x)=0.9 + 0.3 \exp\left(-40 \left( \frac{x}{400} - \frac{1}{3} \right)^{2} \right) + 0.2 \exp\left(-35 \left( \frac{x}{400} - \frac{2}{3} \right)^{2} \right), $$
the flow is subcritical along the whole channel.

\subsubsection{Smooth transition followed by a hydraulic jump in a long domain}
\label{McDo-pseudo2D-SmoothJump}

In this case, the flow is fixed at the inflow and the water height is prescribed at the outflow. We have the following boundary conditions:
\begin{equation*}
  \left\{\begin{array}{l}
     \text{upstream:}\;q=20\text{ m}^3 \text{s}^{-1},\\
     \text{downstream:}\;h=h_{out}.
  \end{array}\right.
\end{equation*}

The channel is initially dry, with a little puddle downstream, \emph{i.e.} initial conditions are:
\begin{equation*}
 h(x,y)= \max(0, h_{out} + z(400,y) - z(x,y))\;\text{and}\;q=0\;\text{m}^3/\text{s}.
\end{equation*}

With the second channel, we define the mean water height by
\begin{multline*}
h_{ex}(x)=0.9 + 0.25 \left(\exp\left(-\frac{x}{40}\right) - 1 \right) +0.25 \exp \left( 15 \left( \frac{x}{400} - \frac{3}{10} \right) \right) \\ \mbox{ for } 0\mbox{ m} \leq  x  \leq 120 \mbox{ m}, 
\end{multline*}
and
$$h_{ex}(x)=\exp\left(-p (x - x^{\star}) \right) \sum_{i=0}^{M} k_{i} \left( \frac{x-x^{\star}}{x^{\star\star}-x^{\star}} \right)^{i} + \phi(x)\quad \mbox{ for } 120\mbox{ m} \leq  x  \leq 200 \mbox{ m},$$
with:\\[0.1cm]
\begin{tabular}{lll}
$\bullet$ $x^{\star}=120$~m,&
$\bullet$ $x^{\star\star}=400$~m,&
$\bullet$ $M=2$,\\
$\bullet$ $k_{0}=-0.183691$,&
$\bullet$ $k_{1}=1.519577$,&
$\bullet$ $k_{2}=-18.234429$,\\
$\bullet$ $p=0.09$,&
$\bullet$ $\phi(x)=1.5 \exp \left(0.16 \left( \dfrac{x}{400} - 1 \right) \right) - 0.3 \exp \left(2 \left( \dfrac{x}{400} - 1 \right) \right)$.\hspace{-5cm}
\end{tabular}\\[0.1cm]
Starting with a subcritical flow, we get a smooth transition to a supercritical zone, and, through a hydraulic jump, the flow becomes subcritical 
in the remaining of the domain.

\section{Transitory solutions}

In section~\ref{sec:steady-state}, we gave steady-state solutions of increasing difficulties.
These solutions can be used to check if the numerical methods are able to keep/catch steady-state flows.
But even if the initial condition differs from the expected steady state,
we do not have information about the transitory behavior.
Thus, in this section, we list transitory solutions  that may improve the validation of the numerical methods.
Moreover, as most of these cases have wet/dry transitions,
one can check the ability of the schemes to capture the evolution of these fronts
 (\emph{e.g.} some methods may fail and give negative water height). At last, we give some periodic transitory solutions in order to
check whether the schemes are numerically diffusive or not.
 
Table~\ref{AnaSolTable-Transitory} lists all transitory solutions available in SWASHES and outlines their main features.

\begin{table}[htbp]
\begin{center}
\begin{tiny}
\setlength{\tabcolsep}{1mm}
\begin{tabular}{|c|c|c|c|c|c|c|c|c|c|c|c|}  \cline{1-11}

\multicolumn{4}{|l|}{\emph{Transitory solutions}}&
\multicolumn{3}{c|}{Slope} &
\multicolumn{4}{c|}{Friction} &
 \multicolumn{1}{c}{} \\[0.1cm] \hline
\multicolumn{1}{|c|}{Type} &\multicolumn{1}{|c|}{Description}&\multicolumn{1}{|c|}{\textsection} & \multicolumn{1}{c|}{Reference} &
                     Null &Const. & Var. & Man.& D.-W. & Other & Null & Comments \\\hline

\multirow{6}{*}{Dam breaks}
&
\begin{tabular}{c}Dam break\\on wet domain without friction\end{tabular}
&\ref{Dam-Wet-NoFriction}
&\cite{Stoker57}
&   X  &  &  & & & &X &
Moving shock. 1D \\

\cline{2-12}

&
\begin{tabular}{c}Dam break\\on dry domain without friction\end{tabular}
&\ref{Dam-Dry-NoFriction}
&\cite{Ritter92}
&   X  &  &   & & & & X &
Wet-dry transition. 1D\\

\cline{2-12}

&\begin{tabular}{c}Dam break\\on dry domain with friction\end{tabular}
&\ref{Dam-Dry-WithFriction}
&\cite{Dressler52}
 &   X   &  &    &      &   X & & &
Wet-dry transition. 1D\\

\hline
\multirow{5}{*}{
\begin{tabular}{c}Oscillations\\without damping\end{tabular}
}
&
\begin{tabular}{c}Planar surface\\in a parabola\end{tabular}
&\ref{Thacker-1D}
&\cite{Thacker81}
  &    &   &  X    & & & & X &
Wet-dry transition. 1D\\

\cline{2-12}

&
\begin{tabular}{c}Radially-symmetrical\\paraboloid\end{tabular}
&\ref{Thacker-2D-axisym}
&\cite{Thacker81}
 &   &    &  X    & & & & X &
Wet-dry transition. 2D\\

\cline{2-12}

&
\begin{tabular}{c}Planar surface\\in a paraboloid\end{tabular}
&\ref{Thacker-2D-planar}
&\cite{Thacker81}
 &    &    & X    & & & &X &
Wet-dry transition. 2D\\

\hline
\begin{tabular}{c}Oscillations\\with damping\end{tabular}
&
\begin{tabular}{c} Planar surface\\in a parabola with friction \end{tabular}
& \ref{Sampson-1D}
& \cite{Sampson06}
&     &   & X                   &    & & X & &  Wet-dry transition. 1D\\
   
\hline
\multicolumn{12}{l}{\emph{Const.: Constant;  Var.: Variable; Man.: Manning; D.-W.: Darcy-Weisbach}}\\

\end{tabular}
\end{tiny}
\end{center}

\hspace*{9em}
{\caption{Analytic solutions for shallow flow equations and their main features --- Transitory cases
\label{AnaSolTable-Transitory}}}

\end{table}


\subsection{Dam breaks}
\label{Dam}

In this section, we are interested in dam break solutions of increasing complexity on a flat topography namely Stoker's,
Ritter's and Dressler's solutions. The analysis of dam break flow is part of dam design and safety analysis: dam breaks can release
an enormous amount of water in a short time. This
 could be a threat to human life and to the infrastructures. To quantify the associated risk, a detailed description of the
 dam break flood wave is required. Research on dam break started more than a century ago. In 1892, Ritter
 was the first who studied the problem, deriving an analytic solution based on the characteristics method
 (all the following solutions are generalizations of his method). 
He gave the solution for a dam break on a dry bed without friction
(in particular, he considered an ideal fluid flow at the wavefront):
 it gives a parabolic water height profile connecting the upstream undisturbed region
 to the wet/dry transition point.
In the 1950's, Dressler (see~\cite{Dressler52}) and Whitham (see~\cite{Whitham55}) derived analytic expressions
 for dam breaks on a dry bed including the effect of bed resistance with Ch\'ezy friction law.
 They both proved that the solution is equal to Ritter's solution behind the wave tip.
 But Dressler neglects the tip region (so his solution gives the location of the tip but not its shape)
 whereas Whitham's approach, by treating the tip thanks to an integral method, is more complete.
Dressler compared these two solutions on experimental data \cite{Dressler54}.
A few years later, Stoker generalized Ritter's solution for a wet bed downstream the dam to avoid
 wet/dry transition.

 Let us mention some other dam break solutions but we do not detail their expressions here.
 Ritter's solution has been generalized to a trapezoidal cross section
 channel thanks to Taylor's series in \cite{Wu86}.
 Dam break flows are also examined for problems in
 hydraulic or coastal engineering for example to discuss the behavior of a strong bore, caused by a tsunami, over a uniformly sloping
 beach. Thus in 1985 Matsutomi gave a solution of a dam break on a uniformly sloping bottom (as mentioned in
 \cite{Matsutomi03}).
 Another contribution is the one of Chanson, who generalized Dressler's and Whitham's dam break solutions to turbulent and laminar flows with
 horizontal and sloping bottom \cite{Chanson05,Chanson06}.


\subsubsection{Dam break on a wet domain without friction}
\label{Dam-Wet-NoFriction}

\noindent

In the shallow water community, Stoker's solution or dam break on a wet domain is a classical case (introduced
 first in \cite[p.~333--341]{Stoker57}).
This is a classical Riemann problem: its analog in compressible gas dynamics is the Sod tube \cite{Sod78} and in blood flow
 dynamics the ideal tourniquet \cite{Delestre11}.
In this section, we consider an ideal dam break on a wet domain, \emph{i.e.} the dam break is instanteneous,
 the bottom is flat and there is no friction. We obtain an analytic solution of this case thanks to the characteristics
 method.\\

The initial condition for this configuration is the following Riemann problem
$$h(x)= \begin{cases}
h_{l} & \mbox{ for }0\mbox{ m}\leq x \leq x_{0} ,\\
h_{r} & \mbox{  for } x_{0} <x\leq L,
\end{cases}$$ with $h_{l} \geq h_{r}$ and $u(x)= 0$~m/s. \\

At time $t\geq 0$, we have a left-going rarefaction wave (or a part of parabola between $x_A(t)$ and $x_B(t)$) that reduces the
 initial depth $h_l$ into $h_m$, and a right-going shock (located in $x_C(t)$) that increases the intial height $h_r$ into
 $h_m$. For each time $t\geq 0$, the analytic solution is given by\begin{equation*}
  h(t,x) =\left\{
  \begin{array}{l}
    h_{l} \\[0.3cm]
    \dfrac{4}{9g}\left(\sqrt{gh_{l}}-\dfrac{x-x_0}{2t}\right)^2  \\[0.5cm]
    \dfrac{c_{m} ^2}{g} \\[0.3cm]
    h_{r}
  \end{array}\right. 
  u(t,x) = \left\{
  \begin{array}{ll}
    0 \mbox{ m/s}  & \text{if }\; x\le {x_A}(t), \\[0.3cm]
    \dfrac{2}{3}\left(\dfrac{x-x_0}{t}+\sqrt{gh_{l}}\right) & \text{if }\; x_A(t)\le x\le x_B(t), \\[0.5cm]
    2\left(\sqrt{gh_{l}}-c_{m}\right)  & \text{if }\; x_B(t)\le x \le x_{C}(t), \\[0.3cm]
    0\mbox{ m/s} & \text{if }\; x_C(t)\le x,
  \end{array}\right.
\end{equation*}
where
\begin{equation*}
  x_A(t)= x_{0}- t\sqrt{g h_{l}}, \  x_B(t)=x_{0}+ t\left(2 \sqrt{g h_{l}} - 3 c_{m}\right)  \mbox{ and }
  x_{C}(t) = x_{0} + t  \dfrac{2c_{m}^2 \left( \sqrt{gh_{l}} - c_{m} \right)}{c_{m}^2 - gh_{r}},
\end{equation*}
with $c_{m}=\sqrt{gh_m}$ solution of $-8gh_{r}{c_m}^2\left(\sqrt{gh_{l}}-{c_m}\right)^2+\left({c_m}^2-gh_{r}\right)^2\left({c_m}^2+gh_{r}\right)=0$.\\

This solution tests whether the code gives the location of the moving shock properly.

In SWASHES, we take the following parameters for the dam: $h_{l} = 0.005$~m, $h_{r} = 0.001$~m,  $x_{0} = 5$~m, $L =10$~m and
 $T=6$~s.


\subsubsection{Dam break on a dry domain without friction}
\label{Dam-Dry-NoFriction}

\noindent
Let us now look at Ritter's solution \cite{Ritter92}:
this is an ideal dam break (with a reservoir of constant height
 $h_{l}$) on a dry domain, \emph{i.e.} as for the Stoker's solution, the dam break is instantaneous, the bottom is flat and
 there is no friction. The initial condition (Riemann problem) is modified and reads:
\begin{equation*}
h(x)= \left\{\begin{array}{ll}
  h_{l} >0 & \mbox{ for } 0\mbox{ m}\leq x \leq x_{0},\\
  h_{r} = 0  \mbox{ m} & \mbox{ for } x_{0}<x\leq L,
\end{array}\right.
\end{equation*}
 and $u(x)=0$~m/s.
At time $t>0$, the free surface is the constant water height ($h_l$) at rest connected to a dry zone ($h_r$) by a parabola. This
 parabola is limited upstream (resp. downstream) by the abscissa $x_A(t)$ (resp. $x_B(t)$). The analytic solution is given by
\begin{equation*}
  h(t,x) =\left\{
  \begin{array}{l}
    h_{l}   \\[0.3cm]
    \dfrac{4}{9g}\left(\sqrt{gh_{l}}-\dfrac{x-x_0}{2t}\right)^2  \\[0.5cm]
    0 \mbox{ m}
  \end{array}\right. 
  u(t,x) = \left\{
  \begin{array}{ll}
    0 \mbox{ m/s}  & \text{if }\; x\le {x_A}(t), \\[0.3cm]
    \dfrac{2}{3}\left(\dfrac{x-x_0}{t}+\sqrt{gh_{l}}\right) & \text{if }\; x_A(t)\le x\le x_B(t), \\[0.5cm]
    0 \mbox{ m/s} & \text{if }\; x_B(t)\le x,
  \end{array}\right.
\end{equation*}
where
\begin{equation*}
  x_A(t)= x_{0}- t\sqrt{g h_{l}} \quad \mbox{and} \quad x_B(t)=x_{0}+ 2t\sqrt{g h_{l}}.
\end{equation*}

This solution shows if the scheme is able to locate and treat correctly the wet/dry transition.
It also emphasizes whether the scheme preserves the positivity of the water height, 
as this property is usually violated near the wetting front.

In SWASHES, we consider the numerical values: $h_{l} = 0.005  \mbox{ m},  x_{0} = 5\mbox{ m}, L =10\mbox{ m}$ and $T=6$~s.


\subsubsection{Dressler's dam break with friction}
\label{Dam-Dry-WithFriction}

In this section, we consider a dam break on a dry domain with a friction term \cite{Dressler52}. 
In the literature we may find several
 approaches for this case. Although it is not complete in the wave tip (behind the wet-dry transition), we present here
 Dressler's approach. Dressler considered Ch\'ezy friction law and used a perturbation method in Ritter's method,
 \emph{i.e.} $u$ and $h$ are expanded as power series in the friction coefficient $C_f=1/C^2$.

\noindent
The initial condition is
\begin{equation*}
h(x)= \left\{\begin{array}{ll}
  h_{l}>0 & \mbox{ for } 0\mbox{ m}\leq x \leq x_{0} ,\\
  h_{r} = 0  \mbox{ m} & \mbox{ for } x_{0} <x\leq L,
\end{array}\right.
\end{equation*}
and $u(x)= 0$~m/s. Dressler's first order developments for the flow resistance give the following corrected water height
 and velocity
\begin{equation}
 \left\{\begin{array}{l}
     h_{co}(x,t)=  \dfrac{1}{g}\left(\dfrac{2}{3}\sqrt{gh_{l}}-\dfrac{x-x_0}{3 t}+\dfrac{g^2}{C^2}\alpha_1 t\right)^2 ,\\[0.5cm]
     u_{co}(x,t)= \dfrac{2\sqrt{gh_{l}}}{3}+\dfrac{2(x-x_0)}{3 t}+\dfrac{g^2}{C^2}\alpha_2 t ,
        \end{array} \right.\label{eq:Dressler-cor}
\end{equation}
where
\begin{equation*}
\alpha_1=
\dfrac{6}{5\left(2-\dfrac{x-x_0}{ t\sqrt{gh_{l}}}\right)}-\dfrac{2}{3}+\dfrac{4\sqrt{3}}{135}\left(2-\dfrac{x-x_0}{ t\sqrt{gh_{l}}}\right)^{3/2}
\end{equation*}
and
\begin{equation*}
\alpha_2=
\dfrac{12}{2-\dfrac{x-x_0}{ t\sqrt{gh_{l}}}}-\dfrac{8}{3}+\dfrac{8\sqrt{3}}{189}\left(2-\dfrac{x-x_0}{ t\sqrt{gh_{l}}}\right)^{3/2}-\dfrac{108}{7\left(2-\dfrac{x-x_0}{t\sqrt{gh_{l}}}\right)^2}.
\end{equation*}
With this approach, four regions are considered: from upstream to downstream, a steady state region ($(h_l,0)$ for $x\leq x_A(t)$),
 a corrected region ($(h_{co},u_{co})$ for $x_A(t)\leq x \leq x_{T}(t)$), the tip region (for $x_{T}(t)\leq x \leq x_B(t)$)
 and the dry region ($(0,0)$ for $x_B(t)\leq x$). In the tip region, friction term is preponderant thus \eqref{eq:Dressler-cor}
 is no more valid. In the corrected region, the velocity increases with $x$. Dressler assumed that at $x_T(t)$ the velocity
 reaches the maximum of $u_{co}$ and that the velocity is constant in space
 in the tip region $u_{tip}(t)=\max_{x\in [x_A(t),x_B(t)]} u_{co}(x,t)$.

The analytic solution is then given by
\begin{equation*}
  h(t,x) =\left\{
  \begin{array}{l}
    h_{l}  \\[0.3cm]
    h_{co}(x,t) \\[0.3cm]
    h_{co}(x,t)  \\[0.3cm]
    0 \mbox{ m}
  \end{array}\right.  \quad
  u(t,x) = \left\{
  \begin{array}{ll}
    0 \mbox{ m/s}  & \text{if }\; x\le {x_A}(t), \\[0.3cm]
    u_{co}(x,t) & \text{if }\; x_A(t)\le x\le x_T(t), \\[0.3cm]
    u_{tip}(t)  & \text{if }\; x_T(t)\le x\le x_B(t), \\[0.3cm]
    0 \mbox{ m/s} & \text{if }\; x_B(t)\le x,
  \end{array}\right.
\end{equation*}

and with
\begin{equation*}
  x_A(t)= x_{0}- t\sqrt{g h_{l}} \quad \mbox{and} \quad x_B(t)=x_{0}+ 2t\sqrt{g h_{l}}.
\end{equation*}

We should remark that with this approach the water height is not modified in the tip zone. This is a limit of Dressler's
 approach. Thus we coded the second order interpolation used in \cite{Valiani99, Valiani02} (not detailed here) and
 recommanded by Valerio Caleffi\footnote{Personal communication}. 
 
 Even if we have no information concerning the shape of
 the wave tip, this case shows if the scheme is able to locate and treat correctly the wet/dry transition.

In SWASHES, we have $h_{l}= 6 \mbox{ m}, x_{0} = 1000 \mbox{ m}$,
 $C=40\mbox{ m}^{1/2}\mbox{/s}$ (Chezy coefficient), $L= 2000\mbox{ m }$ and $T=40$~s.


\subsection{Oscillations}
\label{Thacker}
\vspace{-2pt}

In this section, we are interested in Thacker's and Sampson's solutions. 
These are analytic solutions with
a variable slope (in space) for which the wet/dry transitions are moving. Such moving-boundary solutions are of great interest in communities interested
 in tsunami run-up and ocean flow simulations (see among others \cite{Marche05}, \cite{Kim07} and \cite{Popinet11}). A prime
 motivation for these solutions is to provide tests for numerical techniques and codes in wet/dry transitions on varying
 topographies. The first moving-boundary solutions of Shallow-Water equations
 for a water wave climbing a linearly sloping beach is obtained in~\cite{Carrier58} (thanks to a hodograph transformation).
 Using Shallow-Water equations in Lagrangian form,
Miles and Ball \cite{Miles63} and \cite{Ball63} mentioned exact moving-boundary solutions in a parabolic trough and in a paraboloid of revolution.
 In \cite{Thacker81}, Thacker shows exact moving boundary solutions
 using Eulerian equations. His approach was first to make assumptions about the nature of the motion and then to solve the basin
 shape in which that motion is possible. His solutions are analytic periodic solutions (there is no damping) with Coriolis
 effect. Some of these solutions were generalised by Sampson \emph{et al.} \cite{Sampson06,Sampson07b} by adding damping due to a linear
 friction term.
 
Thacker's solutions described here do not take into account Coriolis effect. In the 1D case, the topography is a parabola and in the 2D case it is a paraboloid.
 These solutions test the ability of schemes to simulate flows with comings and goings and, as the water height is
 periodic in time, the numerical diffusion of the scheme. Similarly to the dam breaks on a dry domain, Thacker's solutions also test
wet/dry transition.

\subsubsection{Planar surface in a parabola without friction}
\label{Thacker-1D}

\noindent
Each solution written by Thacker has two dimensions in space \cite{Thacker81}. The exact solution described here is a simplification to 1D of
 an artificially 2D Thacker's solution. Indeed, for this solution, Thacker considered an infinite channel with a parabolic cross
 section but the velocity has only one nonzero component (orthogonal to the axis of the channel). 
This case provides us with a relevant test in 1D for shallow
 water model because it deals with a sloping bed as well as with wetting and drying. 
 The topography is a parabolic bowl given by
\begin{equation*}
z(x)=h_{0}\left(\dfrac{1}{a^2}\left(x-\dfrac{L}{2}\right)^2-1\right),
\end{equation*}
 and the initial condition on the water height is
\begin{equation*}
h(x)=  \left\{\begin{array}{ll}
  - h_0\left(\left(\dfrac{1}{a}\left(x-\dfrac{L}{2}\right)+\dfrac{B}{\sqrt{2g h_0}}\right)^2-1\right)  &
 \mbox{ for } x_{1}(0)\leq x \leq x_{2}(0),\\
  0  \mbox{ m} & \mbox{ otherwise,}
  \end{array}\right.
\end{equation*}
with $B = \sqrt{2g h_0}/(2a)$
and for the velocity $u(x)= 0$~m/s. 
Thacker's solution is a periodic solution
 (without friction) and the free surface remains planar in time.
The analytic solution is 
\begin{equation*}
  h(t,x) =\left\{
  \begin{array}{ll}
     - h_0\left(\left(\dfrac{1}{a}\left(x-\dfrac{L}{2}\right)+\dfrac{B}{\sqrt{2g h_0}}\cos \left(\dfrac{ \sqrt{2g h_{0}}}{a}t\right)\right)^2-1\right) & \\
     & \hspace{-3cm} \mbox{ for } x_{1}(t)\leq x \leq x_{2}(t), \\[0.5cm]
    0 \mbox{ m}& \hspace{-3cm} \mbox{ otherwise, }
  \end{array}\right.
\end{equation*}
\begin{equation*}  
  u(t,x) = \left\{
  \begin{array}{ll}
    B \sin \left(\dfrac{ \sqrt{2g h_{0}}}{a}t\right) &  \mbox{ for } x_{1}(t)\leq x \leq x_{2}(t),\\[0.5cm]
    0 \mbox{ m/s} & \mbox{ otherwise. }
  \end{array}\right.
\end{equation*}
where $x_1(t)$ and $x_2(t)$ are the locations of wet/dry interfaces at time $t$
\begin{equation*}
 \begin{array}{l}
  x_1(t)=-\dfrac{1}{2}\cos\left(\dfrac{\sqrt{2gh_0}}{a}t\right)-a+\dfrac{L}{2},\\
\\
  x_2(t)=-\dfrac{1}{2}\cos\left(\dfrac{\sqrt{2gh_0}}{a}t\right)+a+\dfrac{L}{2}.
 \end{array}
\end{equation*}

In SWASHES, we consider $a=1$~m, $h_{0}=0.5$~m, $L=4$~m and $T=10.0303$~s (5 periods).


\subsubsection{Two dimensional cases}
\label{Thacker-2D}

Several two dimensional exact solutions with moving boundaries were developed by Thacker. Most of them include the Coriolis
 force that we do not consider here (for further information, see \cite{Thacker81}). These solutions are periodic in time with
 moving wet/dry transitions. They provide perfect tests for shallow water as they deal with bed slope and wetting/drying
 with two dimensional effects. Moreover, as the solution is exact without discontinuity, it is very appropriate to verify
 the accuracy of a numerical method.

\paragraph{Radially-symmetrical paraboloid}
\label{Thacker-2D-axisym}

\begin{figure}
\begin{center}
\includegraphics[width = 0.6 \textwidth]{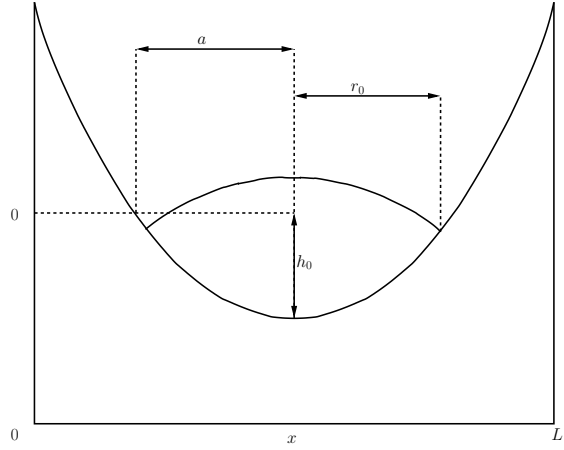}
\caption{Notations for Thacker's axisymmetrical solution}
\label{figThackerAxi}
\end{center}
\end{figure}

The two dimensional case presented here is a radially symmetrical oscillating paraboloid \cite{Thacker81}. The solution is periodic
 without damping (\emph{i.e.} no friction). The topography is a paraboloid of revolution defined by
\begin{equation}
z(r)= -h_{0} \left( 1 -\dfrac{r^{2}}{a^{2}} \right) \label{topoThacker}
\end{equation}
with $r = \sqrt{(x-L/2)^{2}+(y-L/2)^{2}}$ for each $(x,y)$ in $[0;L] \times [0;L]$, where $h_{0}$ is the water depth at
 the central point of the domain for a zero elevation and $a$ is the distance from this central point to the zero elevation of the shoreline.
 The solution is given by:
\begin{equation*}
\left\{\begin{array}{l}
h(r,t) = h_{0} \left( \dfrac{\sqrt{1-A^{2}}}{1-A \cos(\omega t)} - 1 - \dfrac{r^2}{a^2} \left( \dfrac{1-A^2}{\left( 1- A \cos(\omega t) \right)^2} -1 \right) \right) - z(r)\\
u(x,y,t) = \dfrac{1}{1- A \cos(\omega t)} \left( \dfrac{1}{2} \omega \left(x-\dfrac{L}{2}\right) A \sin(\omega t) \right)\\
v(x,y,t) = \dfrac{1}{1- A \cos(\omega t)} \left( \dfrac{1}{2} \omega \left(y-\dfrac{L}{2}\right) A \sin(\omega t) \right)  
\end{array}\right.
\end{equation*}
where the frequency $\omega$ is defined as $\omega = \sqrt{8 g h_{0}}/a$, $r_{0}$ is the distance from the central
 point to the point where the shoreline is initially located and $A = (a^{2}-r_{0}^{2})/(a^{2}+r_{0}^{2})$ (Figure~\ref{figThackerAxi}).\\
The analytic solution at $t=0$~s is taken as initial condition.

In SWASHES, we consider $a=1$~m, $r_{0}=0.8$~m, $h_{0}=0.1$~m, $L=4$~m and $T=3\frac{2\pi}{\omega}$.

\paragraph{Planar surface in a paraboloid}
\label{Thacker-2D-planar}

For this second Thacker's 2D case, the moving shoreline is a circle
and the topography is again given by \eqref{topoThacker}. The free surface has a periodic motion and
 remains planar in time \cite{Thacker81}.
To visualize this case, one can think of a glass with some liquid in rotation inside. 

The exact periodic solution is given by:
\begin{equation*}
 \left\{\begin{array}{l}
h(x,y,t) = \dfrac{\eta h_{0}}{a^2} \left(2 \left(x-\dfrac{L}{2}\right) \cos(\omega t) + 2 \left(y-\dfrac{L}{2}\right) \sin(\omega t) - \eta \right) - z(x,y)\\
u(x,y,t) = -\eta \omega \sin(\omega t)\\
v(x,y,t) = \eta \omega \cos(\omega t)
 \end{array}\right.
\end{equation*}
for each $(x,y)$ in $[0;L] \times [0;L]$, where the frequency $\omega$ is defined as $\omega = \sqrt{2 g h_0}/{a}$ and $\eta$ is a parameter.\\
Here again, the analytic solution at $t=0$~s is taken as initial condition.

In SWASHES, we consider $a=1$~m, $h_{0}=0.1$~m, $\eta=0.5$, $L=4$~m and $T=3\frac{2\pi}{\omega}$.

\subsubsection{Planar surface in a parabola with friction}
\label{Sampson-1D}

Considering a linear friction term (\emph{i.e.} $S_f=\tau u/g$) in system \eqref{SaintVenant1D} (with $R=0$~m/s and $\mu=0~\text{m}^2/\text{s}$) with
Thacker's approach, Sampson \emph{et al.} got moving boundaries solutions with damping  \cite{Sampson06,Sampson07b}. These solutions provide a set of 1D benchmarks for numerical
 techniques in wet/dry transitions on varying topographies (as Thacker's solutions) and with a friction term. One of these solutions
 is presented here. The topography is a parabolic bowl given by
\begin{equation*}
 z(x)=h_0\dfrac{\left(x-\dfrac{L}{2}\right)^2}{a^2},
\end{equation*}
where $h_0$ and $a$ are two parameters and $x \in [0, L]$. The initial free surface is
\begin{equation*}
 (z+h)(x)=\left\{\begin{array}{ll}
                  h_0+\dfrac{a^2 B^2}{8g^2 h_0}\left(\dfrac{\tau^2}{4}-s^2\right) -\dfrac{B^2}{4g}-\dfrac{1}{g}B s \left(x - \dfrac{L}{2}\right) &\\
                  & \hspace{-3cm} \text{for}\; x_1(0)\leq x \leq x_2(0), \\
      0\;\text{m} & \hspace{-3cm} \text{else},
                 \end{array}\right.
\end{equation*}
with $B$ a constant, $s=\sqrt{p^2-\tau^2}/2$ and $p=\sqrt{8gh_0/a^2}$. The free surface remains planar in time
\begin{equation*}
(z+h)(t,x)=\left\{\begin{array}{ll}
        h_0+\dfrac{a^2 B^2 e^{-\tau t}}{8g^2 h_0}\left(-s\tau \sin(2st)+
\left(\dfrac{\tau^2}{4}-s^2\right)\cos(2st)\right) & \\[0.3cm]
\quad-\dfrac{B^2 e^{-\tau t}}{4g}-\dfrac{e^{-\tau t/2}}{g}\left(Bs\cos(st)+\dfrac{\tau B}{2}\sin(st)\right) \left(x - \dfrac{L}{2} \right)&\\ 
& \hspace{-5cm}\text{for}\; x_1(t)\leq x \leq x_2(t), \\[0.4cm]
z(x) &\hspace{-5cm} \text{else}
                  \end{array}\right.
\end{equation*}
and the velocity is given by
\begin{equation*}
u(t,x)=\left\{\begin{array}{ll}
        B e^{-\tau t/2}\sin(s t) & \text{for}\; x_1(t)\leq x \leq x_2(t), \\
0 \;\text{m/s} & \text{else}.
                  \end{array}\right.
\end{equation*}

The wet/dry transitions are located at $x_1(t)$ and $x_2(t)$
\begin{equation*}
\begin{array}{l}
 x_1(t)=\dfrac{a^2 e^{-\tau t/2}}{2gh_0}\left(-Bs \cos(st)-\dfrac{\tau B}{2}\sin(st)\right)-a+\dfrac{L}{2}, \\
\\
 x_2(t)=\dfrac{a^2 e^{-\tau t/2}}{2gh_0}\left(-Bs \cos(st)-\dfrac{\tau B}{2}\sin(st)\right)+a+\dfrac{L}{2}.      
      \end{array}
\end{equation*}

In SWASHES, we consider $a=3,000$~m, $h_0=10$~m, $\tau=0.001\;\text{s}^{-1}$, $B=5$~{m/s},
 $L=10,000$~m and $T=6,000$~s.

\section{The SWASHES software} \label{sec:SWASHES}

In this section, we describe the Shallow Water Analytic Solutions for Hydraulic and Environmental Studies (SWASHES) software.
At the moment, SWASHES includes all the analytic solutions given in this paper.
The source code is freely available to the community through the SWASHES repository
hosted at \url{http://www.univ-orleans.fr/mapmo/soft/SWASHES}.
It is distributed under CeCILL-V2 (GPL-compatible) free software license. 

When running the program, the user must specify in the command 
line the choice of the solution (namely the dimension, the type, the domain and the number of the solution)
as well as the number of cells for the discretization of the analytic solution. The solution is 
computed and can be redirected in a gnuplot-compatible ASCII file.

SWASHES is written in object-oriented ISO C++.
The program is structured as follows: each type of solutions, 
such as \emph{bump}, \emph{dam\_break}, \emph{Thacker}, etc., is written in a specific class.
This structure gives the opportunity to easily implement a new solution,
whether in a class that already exists (for example a new Mac Donald type solution),
or in a new class. 
Each analytic solution is coded with specific parameters (most of them taken from \cite{Delestre10b}). 
In fact, all the parameters are written in the code, except the number of cells.

We claim that such a library can be useful for developers of Shallow-Water codes to evaluate the performances
and properties of their own code, each analytic solution being a potential piece of benchmark.
Depending on the targeted applications and considering the wide range of flow conditions available among the analytic solutions, 
developers may select a subset of the analytic solutions available in SWASHES.
We recommend not to change the values of the parameters to ease the comparison of numerical methods among different codes. However, 
it may be legitimate to modify these parameters to adapt the case to other specific requirements
(friction coefficient, dam break height, rain intensity, etc.). This can be easily done in SWASHES
but, in that case, the code must be renamed to avoid confusions.


\section{Some numerical results: comparison with FullSWOF approximate solutions}
SWASHES was created because we have been developing a software for the resolution
of Shallow-Water equations, namely FullSWOF, and we wanted to validate it against analytic solutions. 
To illustrate the use of SWASHES in a practical case, we first give a short description of the 1D and 2D codes of FullSWOF and then we 
compare the results of FullSWOF with the analytic solutions.
The comparisons between FullSWOF results and the analytic solutions is based on the
relative error in percentage of the water height, using, of course, the 
analytic solution as a reference. This percentage is positive when FullSWOF overestimates the water height and negative when it underestimates it.

\subsection{The FullSWOF program}

FullSWOF (Full Shallow-Water equations for Overland Flow) is an object-oriented C++ code
(free software and GPL-compatible license CeCILL-V2\footnote{http://www.cecill.info/index.en.html}.
Source code available at \url{http://www.univ-orleans.fr/mapmo/soft/FullSWOF/}) 
developed in the framework of the multidisciplinary project METHODE (see \cite{Delestre08,Delestre10b} and \url{http://www.univ-orleans.fr/mapmo/methode/}). We briefly describe here the principles of the numerical methods used in FullSWOF\_2D. The main strategy consists
in a finite volume method on a structured mesh in two space dimensions. Structured meshes have been chosen because
on the one hand digital topographic maps are often provided on such grids, and, on the other hand, it allows to develop
numerical schemes in one space dimension (implemented in FullSWOF\_1D), the extension to dimension two being then straightforward.
Finite volume method ensures by construction the conservation of the water mass, and is coupled with the hydrostatic reconstruction \cite{Audusse04b, Bouchut04} to deal with the topography source term. This reconstruction preserves the positivity of the water height and provides a well-balanced scheme (notion introduced in \cite{Greenberg96}) \emph{i.e.} it preserves at least  hydrostatic equilibrium \eqref{hydeq} (typically puddles and lakes). Several numerical fluxes and second order reconstructions are implemented. Currently, we recommend, based on \cite{Delestre10b}, to use the second order scheme with MUSCL reconstruction \cite{vanLeer79} and HLL flux \cite{Harten83}. 
FullSWOF is structured in order to ease the implementation of new numerical methods.


\subsection{Examples in one dimension}

In this part, we give the results obtained with FullSWOF\_1D for three one-dimensional cases. 
For these examples, we have 500~cells on the FullSWOF\_1D domain
but 2500 for semi-analytic solutions (see Remark~\ref{rem:raff}).
In the first two figures, 
we also plotted the critical height, in order to show directly whether the flow is fluvial or torrential.

\subsubsection{Transcritical flow with shock}
\label{TransFlowShockBench}
\begin{figure}
\begin{center}
  \includegraphics[width = 0.7 \textwidth]{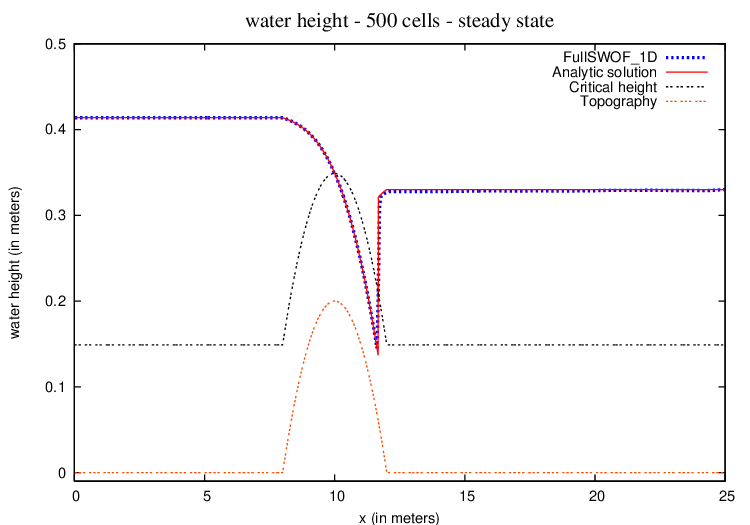}
  \caption{Results of FullSWOF\_1D for a transcritical flow with shock}
  \label{figbump-trans-shock}
\end{center}
\end{figure}

On Figure~\ref{figbump-trans-shock}, we plotted the solution for a transcritical flow with shock 
(Section~\ref{Bump-Transcritical-WithShock-NoFriction}) for a time large enough to attain the steady state, 
namely $t = 100$~s.

Overall, the numerical result of FullSWOF\_1D and the analytic solution are in very good agreement.
By looking carefully at the differences in water height, it appears that the difference is extremely low 
before the bump ($-4\times10^{-6}$~m, \emph{i.e.} -0.001\%). Right at the beginning of the bump, this difference
increases to +0.06\%, but on a single cell. 
Over the bump, differences are alternatively positive and negative which leads, overall, to a good estimate.
A maximum (+1.2\%) is reached exactly at the top of the bump (\emph{i.e.} at the transition from subcritical flow to
supercritical flow). Just before the shock, the difference is only +0.03\%. The largest difference in height (+100\%)
is achieved at the shock and affects a single cell. While for the analytic solution the shock is extremely sharp, 
in FullSWOF\_1D it spans over four cells. After the shock and up to the outlet, the heights computed by FullSWOF\_1D remain
lower than the heights of the analytic solution, going from -1\% (after the shock) to -0.01\% (at the outflow boundary condition).

\subsubsection{Smooth transition and shock in a short domain}
\begin{figure}
\begin{center}
  \includegraphics[width = 0.7 \textwidth]{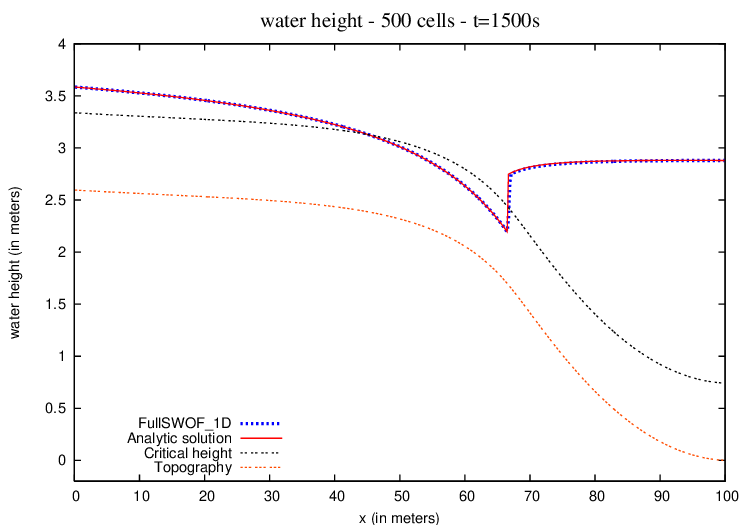}
  \caption{Results of FullSWOF\_1D for the Mac Donald's type solution with 
    a smooth transition and a shock in a short domain, with Manning's friction coefficient}
  \label{figMcDo-100-fluv-Man}
\end{center}
\end{figure}

Figure~\ref{figMcDo-100-fluv-Man} shows the case of a smooth transition and a shock
in a short domain (Section~\ref{McDo-1D-100m-Subcritical}). The final time is $t=150$~s.

FullSWOF\_1D result is very close to the analytic solution.
From $x=0$~m to the subcritical-to-supercritical transition,
FullSWOF\_1D underestimes slightly the water height
and the difference grows smoothly with $x$ (to -0.1\%).
Around $x=45$~m (\emph{i.e.} close to the subcritical-to-supercritical transition), 
the difference remains negative but oscillates, reaching a maximum difference of -0.22\%.
After this transition, FullSWOF\_1D continues to underestimate the water height and this
 difference grows smoothly up to -0.5\% right before the shock.
As for the comparison with the transcritical flow with shock (see Section~\ref{TransFlowShockBench}
), the maximum difference is reached exactly at
the shock ($x=66.6$~m): FullSWOF\_1D overestimates water height by +24\%. Right after the shock,
the overestimation is only +1\% and decreases continuously downstream, reaching 0\% at the outlet.

\subsubsection{Dam break on a dry domain}
\begin{figure}
\begin{center}
  \includegraphics[width = 0.7 \textwidth]{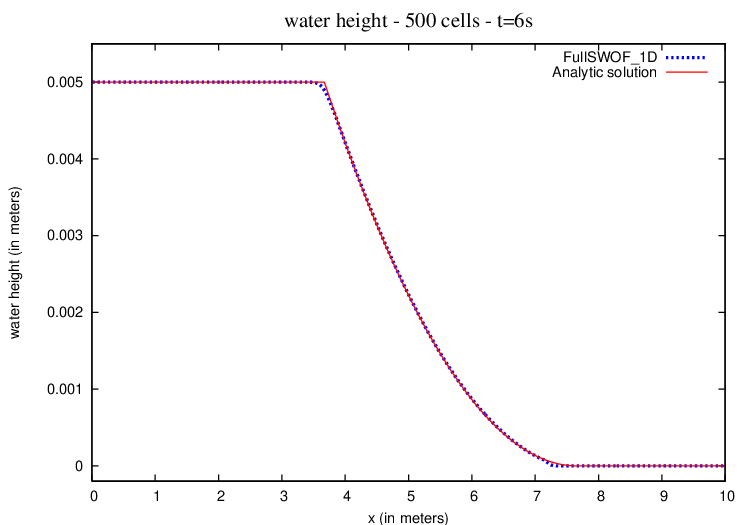}
  \caption{Results of FullSWOF\_1D for a dam break on a dry domain without friction}
  \label{figdam-dry}
\end{center}
\end{figure}

The last figure for the one-dimensional case is Figure~\ref{figdam-dry}, with the solution of 
a dam break on a dry domain (Section~\ref{Dam-Dry-NoFriction}). We chose the final time 
equal to $t=6$~s.
Here too, the numerical result of FullSWOF\_1D matches with the analytic solution well. The height differences 
on the plateau ($h=0.005$ m) are null. On this domain, the water is not flowing yet.
This shows FullSWOF\_1D preserves steady-states at rest properly.
The analytic solution predicts a kink at $x=3.7$~m. In FullSWOF\_1D result, it
does not show a sharp angle but a curve: FullSWOF\_1D slightly underestimates the
water height (by a maximum of -2\%) between $x=3.4$~m and $x=4$~m.
Between $x=4$~m and $x=6.9$~m, FullSWOF\_1D overestimates constantly the water height (by up to +2\%).
On all the domain between $x=6.9$~m and $x=7.6$~m, the water height predicted by FullSWOF\_1D is lower
than the water height computed by the analytic solution (by as much as -100\%).
At this water front, the surface is wet up to $x=7.6$~m according to the analytic solution,
while FullSWOF\_1D predicts a dry surface for $x > 7.3$~m.
It shows the water front predicted by FullSWOF\_1D moves a bit too slowly.
This difference can be due to the degeneracy of the system for vanishing water height.
For $x > 7.6$~m, both the analytic solution and FullSWOF\_1D result give a null water height:
on dry areas, FullSWOF\_1D does not predict positive or negative water heights, and the dry-to-wet
transition does not show any spurious point, contrarily to other codes.


\subsection{Examples in two dimensions}
We now consider the results given by the two-dimensional software FullSWOF\_2D.

\subsubsection{Planar surface in a paraboloid}
\begin{figure}[htbp]
\begin{center}
  \includegraphics[width = 0.7\textwidth]{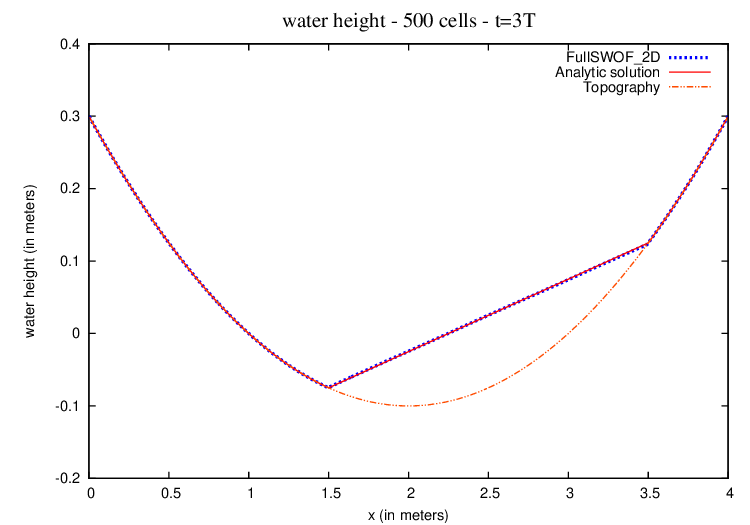}
  \caption{Cross section along~$x$ of FullSWOF\_2D result for Thacker's planar surface in a paraboloid}
  \label{figthacker-2D-plan}
\end{center}
\end{figure}

In the case of Thacker's planar surface in a paraboloid (see Section~\ref{Thacker-2D-planar}),
FullSWOF\_2D was run for three periods on a domain of $[0\text{~m}; 4\text{~m}] \times [0\text{~m}; 4\text{~m}]$
with $500\times 500$~cells. To analyse the performances of FullSWOF\_2D, we consider
a cross section along~$x$ (Figure~\ref{figthacker-2D-plan}).

Overall, FullSWOF\_2D produces a good approximation of the analytic solution:
while the maximum water height is 0.1~m, errors are in the domain
$[-1.55 \times 10^{-3}\text{~m}; +7.65 \times 10^{-4}\text{~m}]$.
At the wet-dry transition located at $x=1.5$~m, FullSWOF\_2D overestimates the water height by about $+6 \times 10^{-4}$~m
and, on a single cell, FullSWOF\_2D predicts water while the analytic solution gives a dry surface.
Close to $x=1.5$~m, the overestimation of the water height is up to +43\%, but decreases quickly and is always
less than +5\% for $x >1.57 $~m. The overestimation persists up to $x=2.65$~m and then becomes an underestimation.
The underestimation tends to grow up to the wet-dry transition at $x=3.5$~m (reaching -5\% at $x=3.38$~m).
Starting exactly at $x=3.5$~m (and on the same cell) both FullSWOF\_2D and the analytic solution predict no water.
However, on the cell just before, FullSWOF\_2D underestimates the water height by -97\% (but at this point, the water height given by the
analytic solution is only $1.6 \times 10^{-3}$~m).

\subsubsection{Mac Donald pseudo-2D solutions}
We consider two Mac Donald pseudo-2D solutions. Since FullSWOF\_2D does not solve the pseudo-2D
Shallow-Water system but the full Shallow-Water system in 2D \eqref{SaintVenant2D}, more significant
differences are expected. In both cases, FullSWOF\_2D was run long enough to reach steady-state.

\paragraph{Supercritical flow in a short domain}

\begin{figure}
\begin{tabular}{cc}
  \includegraphics[height = 7cm, width = 0.51 \textwidth]{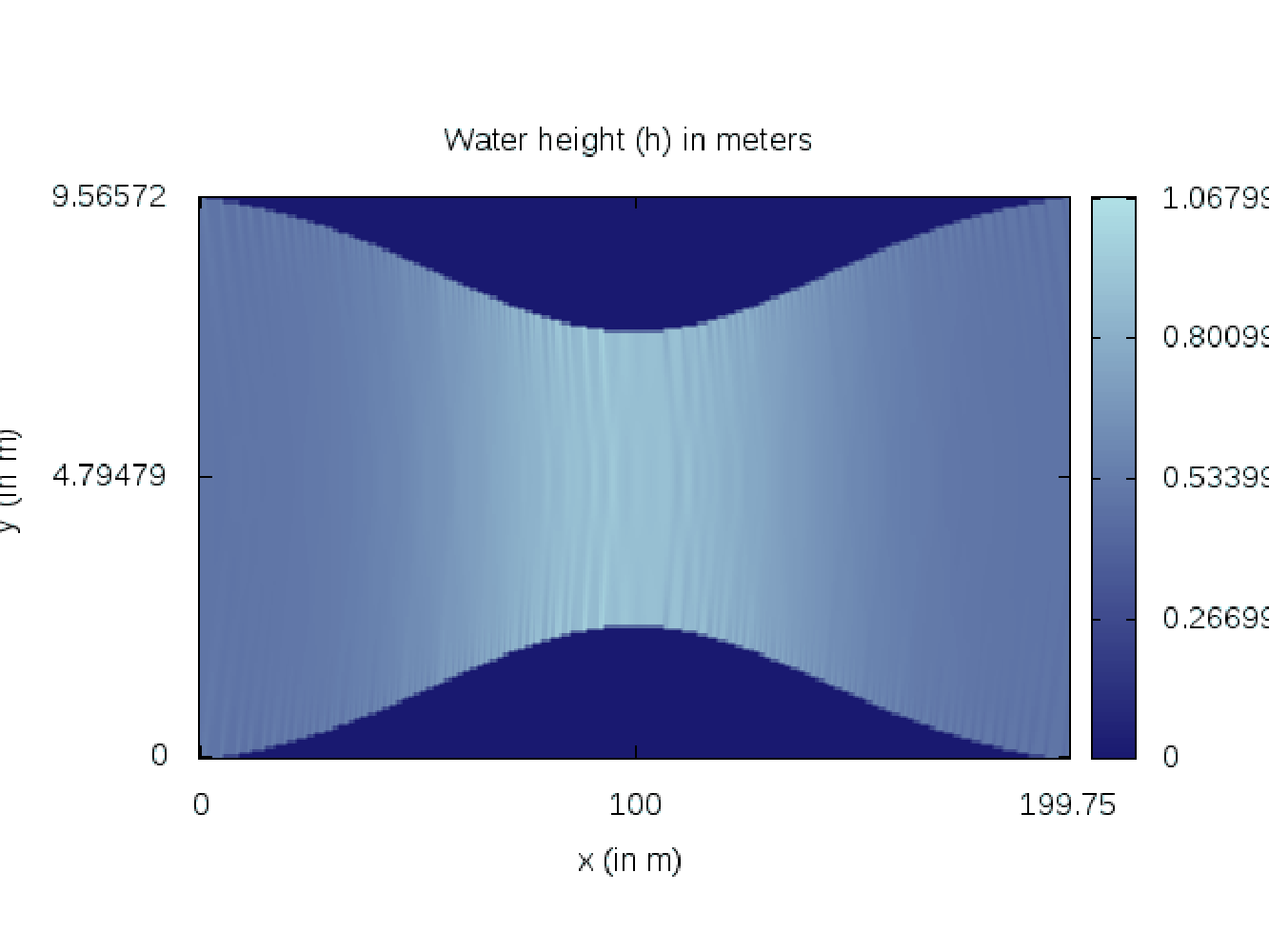} & \\[-6cm]
  &\includegraphics[height = 5cm, width = 0.44 \textwidth]{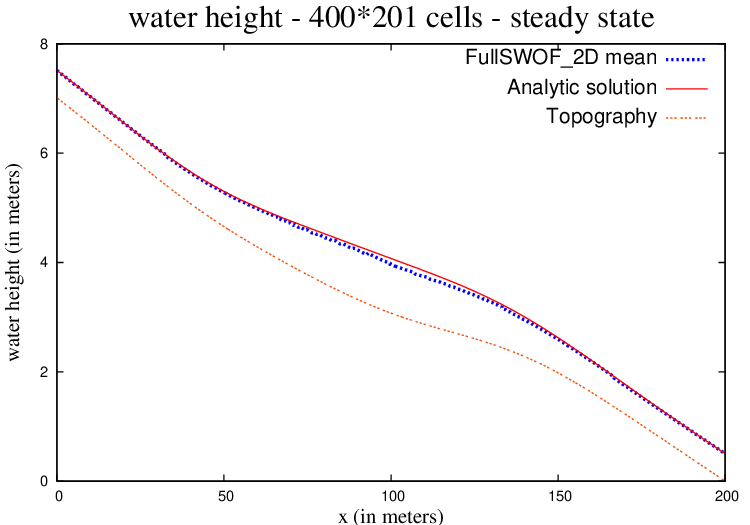}
\end{tabular}
\caption{Results of FullSWOF\_2D for the Mac Donald pseudo-2D supercritical solution in a short domain}
\label{figMacDo-2D-pb10}
\end{figure}

The case of a supercritical flow in a short domain (Section~\ref{McDo-pseudo2D-Supercritical})
is computed with $400 \times 201$ cells on the domain
$[0\text{ m};200\text{ m}] \times [0\text{ m}; 9.589575\text{ m}]$, with the topography $B_1$.
The $y$-averaged result of FullSWOF\_2D differs from the analytic solution mainly around $x = 100$~m
(Figure~\ref{figMacDo-2D-pb10}), with an underestimation of the water height of up
to -0.018~m (-11.85\%). This underestimation occurs on the whole domain and gets closer to zero
near both the upper and lower boundaries.

\begin{figure}
\begin{tabular}{cc}
  \includegraphics[height = 7cm, width = 0.51 \textwidth]{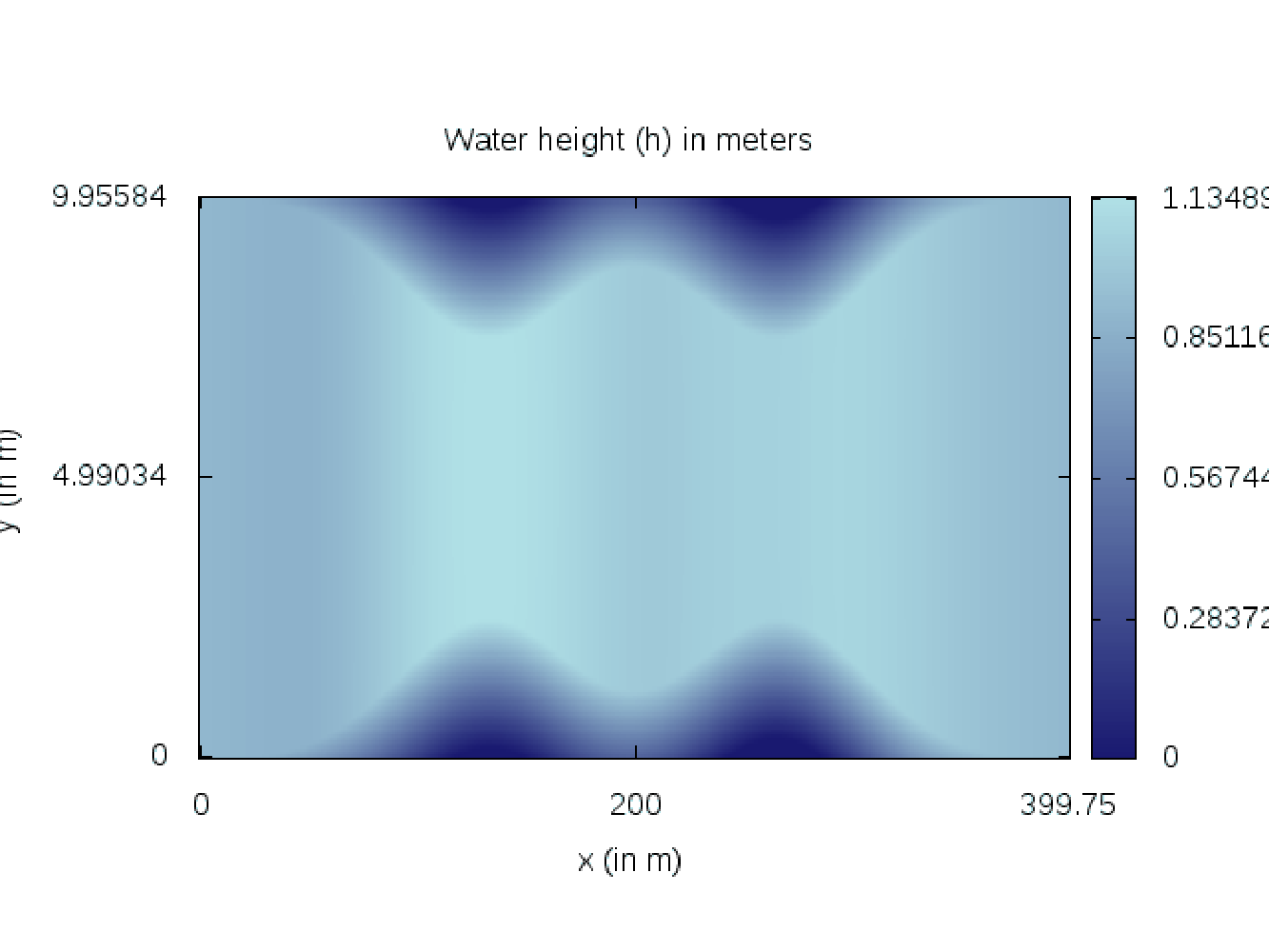} & \\[-6cm]
  &\includegraphics[height = 5cm, width = 0.44 \textwidth]{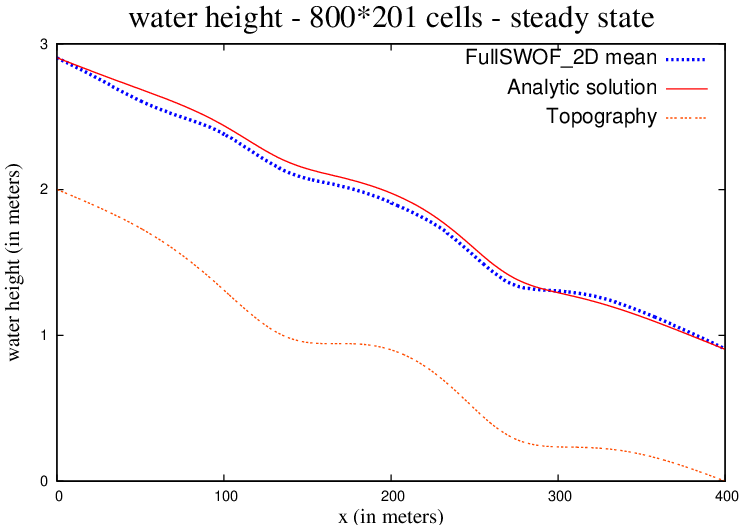}
\end{tabular}
\caption{Result of FullSWOF\_2D for the Mac Donald pseudo-2D subcritical solution in a long domain}
\label{figMacDo-2D-pb13}
\end{figure}

\paragraph{Subcritical flow in a long domain}
The case of a subcritical flow in a long domain (Section~\ref{McDo-pseudo2D-Subcritical2}) 
is computed with $800 \times 201$ cells on the domain
$[0\text{ m};400\text{ m}] \times [0\text{ m};9.98067\text{ m}]$, with the topography $B_2$.
Comparison between the $y$-averaged FullSWOF\_2D result and the analytic solution shows clear
differences (Figure~\ref{figMacDo-2D-pb13}), even if the overall shape of the free surface given
by FullSWOF\_2D matches the analytic solution.
FullSWOF\_2D underestimates the water height on most of the domain.
This difference can be up to -0.088~m (-8.8\%) at $x = 66$~m.
FullSWOF\_2D overestimates water height for $x > 297 $~m and
this overestimation can reach up to +0.04~m (+4.1\%) at $x = 334$~m.

\section*{Acknowledgments}

The authors wish to thank Valerio Caleffi and Anne-C\'eline Boulanger
for their collaboration. This study is part of the ANR METHODE granted by the French National Agency
 for Research {ANR-07-BLAN-0232}.

\bibliographystyle{plain}
\bibliography{SWASHES}

\end{document}